\documentclass[pdflatex,sn-mathphys]{sn-jnl}

\jyear{2021}%

\theoremstyle{thmstyleone}%
\newtheorem{lemma}{Lemma}[section]
\newtheorem{theorem}{Theorem}
\newtheorem{proposition}[theorem]{Proposition}%

\theoremstyle{thmstyletwo}%
\newtheorem{example}{Example}%
\newtheorem{remark}{Remark}%
\newtheorem{corollary}{Corollary}%
\def\d{\mathrm{d}}
\usepackage{subcaption}
\usepackage{bbm}
\usepackage{amsmath}
\usepackage{graphicx}
\def\bbe{{\mathbb E}}
\def\bbq{{\mathbb Q}}
\def\bbp{{\mathbb P}}

\theoremstyle{thmstylethree}%
\newtheorem{definition}{Definition}%
\newcommand{\wh}{\widehat}

\newcommand{\wt}{\widetilde}
\begin{document}

\title[Non-concave optimization]{On the equivalence between Value-at-Risk- and Expected Shortfall-based risk measures in non-concave optimization}

\author[]{\fnm{An} \sur{Chen}}

\author[]{\fnm{Mitja} \sur{Stadje}}

\author*[]{\fnm{Fangyuan} \sur{Zhang}}

\affil[]{\orgdiv{Institute of Insurance Science}, \orgname{University of Ulm}, \orgaddress{\street{Germany}}}

\affil[]{\orgdiv{an.chen@uni-ulm.de (An Chen)}}

\affil[]{\orgdiv{mitja.stadje@uni-ulm.de (Mitja Stadje)}}

\affil[]{\orgdiv{zhang.fangyuan@uni-ulm.de (Fangyuan Zhang)}}


\abstract{We study a non-concave optimization problem in which a financial company maximizes the expected utility of the \emph{surplus} under a risk-based regulatory constraint. 
For this problem, we consider four different prevalent risk constraints (Expected Shortfall, Expected Discounted Shortfall, Value-at-Risk, and Average Value-at-Risk), and investigate their effects on the optimal solution. 
Our main contributions are in obtaining an analytical solution under each of the four risk constraints, in the form of the optimal terminal wealth. 
We show that the four risk constraints lead to the \emph{same} optimal solution, which differs from previous conclusions obtained from the corresponding  concave optimization problem under  a risk constraint. 
Compared with the benchmark (unconstrained) non-concave utility maximization problem, all the four risk constraints effectively and equivalently reduce the set of zero terminal wealth, but do not fully eliminate this set, indicating the success and failure of the respective financial regulations.}

\keywords{Expected Shortfall, Value-at-Risk ,   Average Value-at-Risk,  Non-concave optimization,  Equivalence}



\maketitle

\section{Introduction}\label{sec1}

	Expected utility theory (EUT) is an underlying hypothesis for many models in economics and finance. For instance, EUT serves as a fundamental model for studying the investor's decision under uncertainty \citep{merton1969lifetime,samuelson1975lifetime}, and comparing or pricing financial products \citep{chen2015utility,liu2021indifference}. The influential work by Merton \cite{merton1969lifetime} proposes an optimization problem based on EUT, which became a paradigm model for optimal asset allocation problems in continuous-time settings. Recently, there are many extensions of asset allocation models in order to allow for more practical considerations. For instance, to study the effect of the risk-based regulation in the financial industry, one important stream is to add a risk constraint, e.g. Value-at-Risk and Expected Shortfall, to the optimization problem. In handling the optimization problem under the risk constraints, a dominant and efficient ansatz is the Lagrangian approach from standard convex analysis. For example, the Lagrangian approach has been used to obtain closed-form solutions for the optimal asset allocation problem under a Value-at-Risk (VaR) constraint or a Limited Expected Loss (a.k.a. Expected Discounted Shortfall) constraint \cite{basak2001value}, under a combined VaR-portfolio-insurance constraint \cite{chen2018optimal}, and under multiple-VaR constraints in a multi-period model \cite{chen2018risk}. A recent work uses quantile formulation \citep{he2011portfolio,xu2016note} and the Lagrangian approach to investigate the optimal asset allocation problem under an Average-Value-at-Risk (AVaR) constraint \cite{wei2021risk}.
	
		The other important stream of asset allocation models is to consider a non-concave utility maximization problem. 
	For instance, an $S$-shape utility function is considered to model loss aversion and risk seeking behavior of an agent \citep{tversky1992advances,berkelaar2004optimal,jin2008behavioral}. A non-concave utility function resulting from a non-linear payoff function is modeled to study the managerial compensation problem \citep{carpenter2000does,basak2007optimal}, the participating contract problem \citep{chen2019constrained,nguyen2020nonconcave}, or the general option-type payoff problem \citep{de2011portfolio}. To address the non-concavity, some works first try to recover convexity by a decomposition method \citep{jin2008behavioral,de2011portfolio} or by concavification \citep{carpenter2000does,reichlin2011non}.\footnote{\cite{jin2008behavioral} and \cite{de2011portfolio} decompose the terminal wealth into two random variables $X$ and $Y$ such that the non-concave function is concave with respect to $X$ and convex with respect to $Y$. 
	Then, the key challenge is to determine the proper set on which the optimal solution is defined, which is not intuitive. \cite{carpenter2000does} uses  concavification to solve the unconstrained non-concave optimization problem. Concavification means that one first establishes the concave envelope, e.g. the smallest concave function that dominates the original non-concave function, and then applies standard convex analysis. The optimal solution for the problem with the concave envelope is the same as the optimal solution for the original unconstrained non-concave optimization problem. 
	A rigorous proof for this technique can be found for instance in \cite{reichlin2013utility}.} These methods are inspiring but are also complex to apply. More importantly, if we consider risk constraints, recovering convexity is not sufficient to solve the non-concave optimization problem. 
	
	In this work, we study non-concave optimization under risk constraints, which is related to both extensions above. We consider an optimal asset allocation problem for a managerial board of a financial company, which maximizes the expected utility arising from the \emph{surplus}, i.e., the positive difference between the asset and the liability. As the surplus is a non-linear function of the terminal wealth, the resulting expected utility maximization problem is \emph{non-concave}. In addition, we consider four risk constraints (Expected Shortfall, Expected Discounted Shortfall, Value-at-Risk, Average Value-at-Risk) in the optimization problem. These risk constraints are well known and are discussed widely in both the literature and the industry. Particularly, \emph{concave} utility maximization under these risk constraints have been extensively studied \citep{basak2001value,gandy2005portfolio,cuoco2008optimal,wei2021risk}. However, the closed-form solutions for the non-concave utility maximization problem under these risk constraints are not available in the literature except under a VaR constraint.\footnote{One can find the optimal solution for a non-concave maximization problem arising from a non-linear payoff function under a VaR constraint in \cite{nguyen2020nonconcave}, and the optimal solution with an S-shape utility function under a VaR constraint in \cite{dong2020optimal}.} Consequently, a thoroughly comparative analysis of the effect of these risk constraints on the non-concave utility maximization is missing. Our closed-form solutions for the non-concave utility maximization under these risk constraints fill this gap.\footnote{We remark that there are some discussions about the non-concave optimization under an AVaR constraint but these discussions are restricted to special cases. \cite{de2011portfolio} studies a non-linear payoff function under an AVaR constraint, but the AVaR constraint only applies where the utility is not concave, while we do not have such a restriction. \cite{armstrong2019risk} discusses an S-shape utility function under an AVaR constraint assuming that a loss can be arbitrarily large and that the terminal wealth can be negative. In this case, an AVaR constraint is never efficient in constraining the risk-seeking behavior. On the contrary, we show that an AVaR constraint is effective if the loss is bounded from below and the terminal wealth is non-negative.}
	
	We use the Lagrangian approach to solve the non-concave optimization problems under these constraints.\footnote{Note that the Lagrangian approach can be used to solve a non-concave maximization problem without risk constraints as well, see \cite{berkelaar2004optimal} for the case with an S-shape utility function, and \cite{basak2007optimal,chen2019constrained} for the case with a non-concave utility arising from the non-linear payoff.}  Due to the non-concave utility function and the existence of a risk constraint, the constructed Lagrangian for solving the optimization problem is discontinuous and highly non-concave. By decomposing the non-concave Lagrangian into local concave or affine functions on disjoint sets, we obtain the global maximum by comparing the local maximums of the piece-wise Lagrangian, which reduces to finding the zero roots of a sequence of conjugate functions. However, it is non-trivial to determine the global maximum as the conjugate functions do not necessarily have zero roots, which implies that the existence of the optimal solution is \emph{not} fully ensured. Note that in the concave optimization under these risk constraints, the existence of the optimal solution is guaranteed by the existence of the Lagrangian multipliers. However, in the non-concave optimization under these risk constraints, the optimal solution exists if the Lagrangian multipliers as well as the zero roots of the conjugate functions both exist, which is more challenging to prove. In this paper, we investigate the conditions under which the zero roots of the conjugate functions exist (Lemma \ref{tang}), and obtain the optimal solutions for the non-concave utility maximization problems under these risk constraints (Section 4).  
	
	Compared to the counterpart concave maximization problems, Expected Shortfall and Expected Discounted Shortfall bring more challenges in solving the non-concave utility maximization problem. Note that the Lagrangian multipliers for the risk constraints (Expected Shortfall or an Expected Discounted Shortfall) are explicitly involved in the optimal solution for the concave utility maximization problem. Therefore, it is straightforward to find an optimal binding solution for these problems.\footnote{A non-trivial optimal solution normally satisfies the given risk constraint with an equality, which is also called a binding solution.} Contrary to our intuition, the Lagrangian multiplier for those risk constraints (Expected Shortfall and Expected Discounted Shortfall) is no longer involved in the solution, and thus finding an optimal binding solution is non-trivial. Consequently, determining the corresponding global maximum is a substantially more complex problem compared to the case under a VaR constraint in a non-concave model (Appendix B).
	
	The optimization problem under an AVaR constraint is more challenging regardless of the utility function being concave or not. The reason is that an AVaR constraint is formulated based on the quantile function of the terminal wealth, while the optimization problem per se is formulated based on the terminal wealth. We deal with this problem in two steps. In the first step we transform an AVaR constraint to an equivalent set of Expected Shortfall constraints (Lemma \ref{lemma:AVaR}), which transforms the quantile-based constraint to the terminal wealth-based constraint. Then, we show that the optimal solution for a non-concave utility maximization problem under an AVaR constraint is equivalent to the optimal solution under the equivalent Expected Shortfall constraint (Proposition \ref{proposition:AVaR}).\footnote{Recently, the optimal solution for a \emph{concave} utility maximization problem under an AVaR constraint is provided by Wei \cite{wei2021risk}. He first formulates an equivalent optimization problem based on the quantile function of the terminal wealth via quantile formulation \citep{he2011portfolio,xu2016note}, and then applies the Lagrangian approach to solve the optimization problem under an AVaR constraint.}
	
	Based on the closed-form solutions, we compare these different risk constraints on the non-concave utility maximization to evaluate the efficiency of the risk-based regulation, assuming that the financial company is \emph{surplus-driven}. Note that in the counterpart \emph{concave} utility maximization problems, the regulatory risk constraints discussed, mainly affect the tail distribution of the optimal wealth. Moreover, it has been shown that the Expected Discounted Shortfall constraint leads to a higher terminal wealth in the tail than the unconstrained case, while the other three constraints lead to a lower terminal wealth than the unconstrained case in a concave model \citep{basak2001value,wei2021risk}. However in a non-concave model, the four risk constraints, though defined differently, all lead to the \emph{same optimal wealth} in most relevant situations. 
In fact, it has been noticed in the literature that the company with a non-concave utility function shows a gambling behavior, by incurring a non-negligible set of zero terminal wealth or the maximum possible loss, see \cite{carpenter2000does,berkelaar2004optimal,jin2008behavioral}. Our results show that all the four popular risk constraints \emph{cannot} eliminate the set of zero terminal wealth, but can reduce the set of zero terminal wealth compared with the unconstrained case. Therefore, once considering the potential surplus-driven characteristic of the financial company, the four risk constraints are \emph{as good as} or \emph{as bad as} each other. Although the optimization problem is formulated to study the optimal asset allocation problem of a surplus-driven financial company under the risk-based regulation, the result in this paper can be easily extended to the area of e.g. the managerial compensation and the participating contracts problems by slightly adjusting the non-linear payoff functions.

The paper is arranged as follows. In Section 2, we define the four risk constraints formally and formulate the non-concave utility maximization problems under these constraints in a complete financial market. In Section 3, we briefly summarize the existing solutions for the counterpart concave utility maximization problems under these constraints. In Section 4, we provide the optimal solutions for the non-concave utility maximization problems under the risk constraints. In Section 5, we discuss the equivalence among these optimal solutions, and illustrate the equivalence by numerical examples. The last section concludes shortly. All technical proofs are put in appendices. Moreover, we provide a general procedure for solving the non-concave utility maximization problems under risk constraints in Appendix A.
\section{The model}\label{sec2}

	\subsection{The financial market}
	We assume a \emph{complete} financial market without transaction costs in continuous time that contains one traded risk-free asset $S_0$ (the bank account) and $m$ traded risky assets denoted by the stochastic processes $S=(S_1,\cdots, S_m)^{'}$.\footnote{Here $S$ is an $m$-dimensional vector and $'$ denotes the transposed sign.} We fix a filtered probability space $(\Omega,\mathcal{F}=(\mathcal{F}_t)_{t\in[0,T]},\bbq)$, $T<\infty$. The unique local martingale measure is denoted by $\bbq$. The state price density process is defined by $\xi_T:=\frac{S_0(0)\d \bbq}{S_0(T)\d \bbp}$.\footnote{The state price density $\xi_T$ is defined in this way such that the discounted asset process is a local martingale under the risk-neutral probability measure $\bbq$.} Throughout the paper, we assume that $\xi_T$ is atom-less.  As the value of $\xi_T$ is high in a recession and low in a prosperous time, $\xi_T$ has the property of directly reflecting the overall state of the economy. Therefore, the functional relationship between the optimal wealth and $\xi_T$ will be used as an interpretation for some of the results.

	The financial institution endowed with an initial capital $x_0$ chooses an investment strategy that we describe by $\pi_i(t)$, the units of $i$th risky asset in the portfolio at time $t$. The strategy is self-financing (i.e., no intermediate income). We assume that $\pi(t)=(\pi_{1}(t),\cdots,\pi_{m}(t))$ is adaptive with respect to the filtration $\mathcal{F}=(\mathcal{F}_t)_{t\in[0,T]}$. The wealth process related to a strategy $\pi(t)$ starting with the initial wealth $x_0$ is then given by
	\begin{align}
	X_t^{\pi}=X_0+\sum_{i=0}^{m}\int_{0}^{t}\pi_i(s)\d S_i(s)=X_0+\int_{0}^{t}\pi(s)\d S(s),\quad X_0=x_0>0.
	\end{align}
	
		In addition, the set of attainable terminal wealth is defined by
	$$ \mathcal{X}:=\{X_T^{\pi}\,\text{is $\mathcal{F}_T$-measurable, replicable, non-negative and}\,\bbe[\xi_TX_T^{\pi}]= x_0 \}.\footnote{$\bbe[]$ denotes that the expectation is taken under the physical probability $\bbp$.}$$
	
		\noindent Note that in a complete financial market, it is sufficient to determine the optimal terminal wealth  from the set of attainable wealth, and the corresponding investment strategy can be obtained by the martingale approach.\footnote{After determining the optimal terminal wealth, one can calculate the optimal investment strategy for a CRRA investor through a standard procedure via It\^o's Lemma in a Black-Scholes market, see the examples in \cite{basak2001value,carpenter2000does,basak2007optimal,chen2018optimal}.} Hence, from now on we omit the dependence of $X_T$ on $\pi$ and focus on finding the optimal terminal wealth. Without ambiguity, we use $\xi$ instead of $\xi_T$ to denote the state price density at time $T$.
		
			\subsection{The model setup with an additional regulatory constraint}
	In this section, we formulate the non-concave utility maximization problem under a risk constraint. We consider four risk constraints, i.e., Expected Shortfall (ESP), Expected Discounted Shortfall (ESQ), Value-at-Risk (VaR) and Average Value-at-Risk (AVaR), which are all frequently discussed and applied in the literature and in practice.\footnote{Expected Shortfall takes the expectation under the physical probability measure $\bbp$, and hence is abbreviated as ESP, while Expected Discounted Shortfall takes the expectation under the risk-neutral measure $\bbq$, and hence is abbreviated as ESQ.} We first introduce the definition of the risk constraints, and then we formulate the optimization problem under each of them, respectively.  
	
		\subsubsection{The definition of the risk constraints}
	In order to ease the notation, we always use $L$ to denote a given regulatory threshold in the different risk constraints. However, the economic interpretations of $L$ differ in each case. 
	\paragraph{Expected Shortfall and Expected Discounted Shortfall}
	Given a regulatory threshold $L$, Expected Shortfall (ESP) is defined by
	\begin{equation}\label{def:es}
	ESP:=\bbe[(L-X_T)^+]\leq\epsilon_\bbp,
	\end{equation}
	where $(L-X_T)^+:=\max(L-X_T,0)$. Defined in this way, the risk constraint ESP \eqref{def:es} aims to control the expected shortfall of the portfolio, where the shortfall occurs if the terminal wealth is below the regulatory threshold $L$.
	
	Similarly, Expected Discounted Shortfall under the same regulatory threshold $L$ is defined by
	\begin{equation}\label{def:esq}
	ESQ:=\bbe^{\bbq}[e^{-rT}(L-X_T)^+]=\bbe[\xi(L-X_T)^+]\leq\epsilon_\bbq.
	\end{equation}
	Obviously, definition \eqref{def:esq} controls the discounted expected shortfall of the portfolio under the risk-neutral measure $\bbq$. We inherit the name ``Expected Discounted Shortfall'' from \cite{shi2012short}. Note that the same risk constraint \eqref{def:esq} is called Limited Expected Loss in \cite{basak2001value}.
	
		\paragraph{Value-at-Risk and Average Value-at-Risk}Value-at-Risk (VaR) is a quantile-based risk measure. For a terminal wealth $X_T$, the VaR constraint at probability $\alpha$ ($0<\alpha< 1$), $VaR_{X_T}(\alpha)$, is defined by
\begin{equation}\label{def:quantile}
		VaR_{X_T}(\alpha):=\sup\{x\vert\bbp(X_T< x)\leq\alpha\}\geq L.
		\end{equation}
		Thus, the VaR constraint requires the upper $\alpha$-quantile of the terminal wealth to be above a given threshold $L$. Note that the VaR constraint, $VaR_{X_T}(\alpha)\geq L$, is equivalent to 
		\begin{equation}\label{def:var}
			\bbp(X_T\geq L)\geq 1-\alpha,
			\end{equation}
			which is consistent with the literature \citep{basak2001value,nguyen2020nonconcave,dong2020optimal}.
		The economical interpretation behind (\ref{def:var}) is that by satisfying the VaR constraint, the terminal wealth of the financial company is above the regulatory threshold $L$ with at least probability $1-\alpha$ (or the terminal wealth drops below the regulatory threshold $L$ with at most probability $\alpha$).\\ 
	
		\noindent Last but not the least, an AVaR constraint is defined by
	\begin{equation}\label{def:AVaR}
			AVaR_{X_T}(\alpha):=\frac{1}{\alpha}\int_{0}^{\alpha} VaR_{X_T}(\beta)\d\beta\geq L^{AVaR},
			\end{equation}	
			where $L^{AVaR}$ is a given regulatory parameter.\footnote{Note that the average quantile in the tail is naturally lower than the $\alpha$-quantile by the definition. Therefore, we use a different notation $L^{AVaR}$ to allow for free choices of the parameters.} The VaR constraint aims to control the $\alpha$-quantile of the portfolio, while Average Value-at-Risk (AVaR) aims to control the average quantile of the portfolio in the tail.
			
					These risk constraints are widely discussed in the literature. However, the same risk constraint can have a different name in other references. In order to avoid ambiguity, we give a remark on the definitions of these risk constraints.
					
						\begin{remark}
			\begin{enumerate}[i)]
				\normalfont
				\item VaR and AVaR are sometimes defined with a minus sign in other works, e.g., \cite{rockafellar2000optimization,acerbi2001expected}. In those contexts, a loss per se is modeled by a (negative) random variable, and the corresponding VaR and AVaR are interpreted as the minimum amount of capital that is needed to make the portfolio acceptable. But in our context, we only focus on the (non-negative) terminal wealth. Consequently, VaR and AVaR are the positive quantile values of the terminal wealth.  Other than the minus sign, the mathematical definitions of VaR and AVaR in our setting are consistent with the literature.
				\item We use the name AVaR to be consistent with \cite{follmer2016stochastic}. It is also called conditional Value-at-Risk \citep{rockafellar2000optimization,rockafellar2002conditional} or Expected Shortfall in other references, e.g., \cite{acerbi2001expected,acerbi2002expected,acerbi2002coherence,embrechts2013model,embrechts2014academic,armstrong2019risk}. Despite the various names, they essentially measure \emph{the average quantile of the portfolio in the tail}. Moreover, an AVaR constraint can be equivalently represented in the form of an Expected Shortfall constraint under specific conditions, see Lemma \ref{lemma:AVaR}. As we also cover the Expected Shortfall (ESP) constraint, we use the name AVaR to distinguish it from the ESP constraint.
					\end{enumerate}
		\end{remark}
		
			\subsubsection{The optimization problem under a risk constraint}	
	We consider a surplus-driven financial institution operating in $[0,T]$, $T<\infty$. At time 0, the company receives an initial contribution, $E_0$, from equity holders and an amount, $D_0$, from debt holders. Consequently, the initial asset value of the company is given by $x_0=E_0+D_0$.
	
		We assume that the benefits of the debt holders are paid out as a lump sum at time $T$. The benefits can be represented as follows: if there is no default at maturity, they correspond to the initial contributions of the debt holders accumulated with a (nominal) rate of return, i.e., $D_T= D_0 e^{\int_0^T g_s ds}$, with $D_T \geq D_0e^{\int_{0}^{T}r_s\d s}$, where $r_s$ is the risk-free rate.\footnote{Note that if $D_T< D_0e^{\int_{0}^{T}r_s\d s}$, the debt holders are better off by investing the money $D_0$ fully in the risk-free asset.} In case of default, i.e., if $X_T \leq D_T$, $X_T$ is paid out at maturity since the company has limited liability. More compactly, the debt holders' terminal payoff is represented as
	$$\varphi_L(X_T)
	=\min (D_T,X_T).$$
	The surplus function of a financial institution, which can also be considered as the payoff to the equity holders, is the residual of the wealth, i.e.,  $$\varphi_E(X_T)=X_T-\varphi_L(X_T)=\max(X_T-D_T,0)=:(X_T-D_T)^+.$$
	
	The surplus-driven company invests the total proceedings $x_0$ in a diversified portfolio of the risky and the risk-free assets as defined in Section 2.1. Further, the utility function of the managerial board is denoted by $U$. We assume that the utility function is defined on the non-negative real line, strictly increasing, strictly concave, bounded from below,\footnote{Unbounded utility functions will be discussed later.} continuously differentiable, and satisfying the usual Inada and asymptotic elasticity (AE) conditions,
	\begin{align*} \label{inada}
	&\textbf{Inada:} \quad U^{'}(0)=\lim_{x\to 0}U^{'}(x)=\infty,\quad U^{'}(\infty)=\lim_{x\to\infty}U^{'}(x)=0,\\
	&\textbf{AE:}\quad\lim_{x\to\infty}\sup\frac{xU^{'}(x)}{U(x)}<1,
	\end{align*}
	\noindent where $U^{'}(x)$ denotes the first derivative of the utility function or the marginal utility. 
	
		\begin{figure}[ht!]
		\centering
		\includegraphics[width=0.5\textwidth]{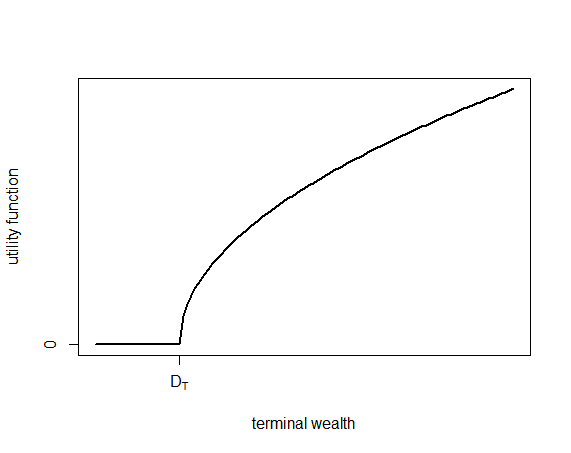}%
		\caption{An illustration of the non-concave utility function of the surplus. As the surplus is zero if the terminal wealth is below the debt level, the utility of the surplus is zero if the terminal wealth is below the debt level.}
		\label{uequ}
	\end{figure}
Figure \ref{uequ} is an illustrative example for the utility function of a surplus-driven company, where the utility level equals zero if the wealth is below the debt level. Although the utility function $U$ is strictly concave, the utility of the surplus $U((X_T-D_T)^+)$ is not globally concave due to the non-linear function of the surplus. Thus, the consequential target function in the optimization problem is not concave. The optimization problem is formulated to maximize the expected utility of the \emph{surplus} by choosing the optimal terminal wealth from the set of the attainable wealth (Section 2.1) while satisfying a risk constraint (ESQ, ESP, VaR or AVaR). We formulate the optimization problems as followings,
	\paragraph{Problem ESQ}
		\normalfont
		\begin{equation*}
		\max\limits_{X_T\in\mathcal{X}}\bbe[U((X_T-D_T)^+)],\quad\mbox{subject to}\quad\bbe[\xi(L-X_T)^+]\leq\epsilon_\bbq,\quad \bbe[X_T\xi]\leq x_0,
		\end{equation*}

\paragraph{Problem ESP}
	\begin{equation*}
	\max\limits_{X_T\in\mathcal{X}}\bbe[U((X_T-D_T)^+)],\quad\mbox{subject to}\quad\bbe[(L-X_T)^+]\leq\epsilon_\bbp,\quad \bbe[X_T\xi]\leq x_0,
	\end{equation*}

\paragraph{Problem VaR}
		\begin{equation*}
		\max\limits_{X_T\in\mathcal{X}}\bbe[U((X_T-D_T)^+)],\quad\mbox{subject to}\quad\bbp(X_T\geq L)\geq 1-\alpha,\quad \bbe[X_T\xi]\leq x_0,
		\end{equation*}
	
\paragraph{Problem AVaR}
	\begin{equation*}
	\max\limits_{X_T\in\mathcal{X}}\bbe[U((X_T-D_T)^+)],\quad\mbox{subject to}\quad\frac{1}{\alpha}\int_{0}^{\alpha}VaR_{X_T}(\beta)\d\beta\geq L^{AVaR},\quad \bbe[X_T\xi]\leq x_0.
	\end{equation*}
	
		In order to better understand the distinct effect of a risk constraint on a non-concave optimization problem, we first review the optimal solution for i) a \emph{non-concave} utility maximization problem \emph{without} a risk constraint; ii) a \emph{concave} optimization problem under a risk constraint. The readers who are familiar with these optimization problems and their solutions can jump to Section 4 directly.
		
			\section{Review on the existing solutions for (non-) concave optimization problems with(-out) constraints}
			
				\subsection{Non-concave optimization \emph{without} a risk constraint}
	If the surplus-driven company does not consider risk constraints, our optimization problem reduces to a non-concave optimization problem without constraints. We refer to this problem as the benchmark problem: 
\paragraph{Benchmark problem}
	\begin{equation*}
		\max\limits_{X_T\in\mathcal{X}}\bbe[U((X_T-D_T)^+)],
		\quad\mbox{\normalfont subject to}\quad \bbe[X_T\xi]\leq x_0.
		\end{equation*}
The benchmark problem can be solved with a similar procedure by concavification \citep{carpenter2000does,reichlin2013utility}, or by the Lagrangian approach \citep{berkelaar2004optimal}. Here, we give the optimal solution without proof.

\begin{figure}[ht]
		\centering
		\includegraphics[width=0.7\textwidth]{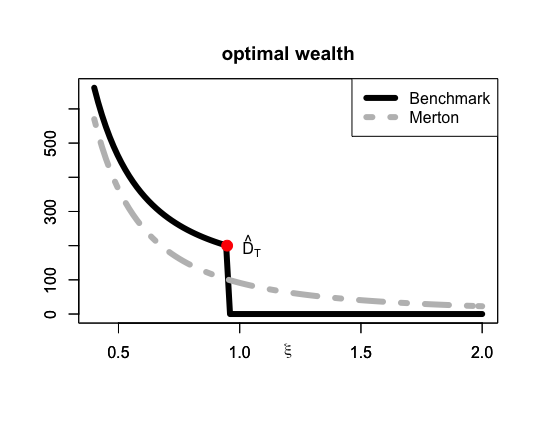}%
		\caption{The optimal terminal wealth of a financial institution without constraints. This figure plots an example for the optimal terminal wealth of a financial institution without constraints. We plot the case of a surplus-driven company as well as the case in which a company maximizes the expected utility of the total asset. The basic parameters for this example are chosen according to Table \ref{table:basic} in Section 5.}	\label{poffB}
	\end{figure}

From now on, let $(U^{'})^{-1}$ denote the inverse function of $U^{'}$ (marginal utility). We introduce the indicator function 
\[
\mathbbm{1}_{\mathcal{A}}=\begin{cases}1,&\quad\mbox{if}\quad \mathcal{A}\quad \mbox{occurs};\\
0,&\quad \mbox{otherwise}.
\end{cases}\]
The optimal solution for the benchmark problem, $X_T^B$, is given by the following proposition.
\begin{proposition}\label{proposition:benchmark}
\begin{enumerate}[i)]
    \item If the utility function is bounded from below, i.e., $U(0)>-\infty$,
    \begin{equation}
        X_T^B=[(U^{'})^{-1}(\lambda_B\xi)+D_T]\mathbbm{1}_{\xi<\xi_B},
    \end{equation}
    where $\lambda_B$ is the Lagrangian multiplier of the budget constraint obtained by solving $\bbe[\xi X_T^B]=x_0$, and $\xi_B$ is the critical value of the state price density defined by $\xi_B=U^{'}(\wh D_T-D_T)/\lambda_B$. $\wh D_T$ is the unique tangent point of function $U(X_T-D_T)$ defined by Lemma A.1. 
    \item If the utility function is not bounded from below, i.e., $U(0)=-\infty$, 
     \begin{equation}
        X_T^B=(U^{'})^{-1}(\lambda_B\xi)+D_T.
    \end{equation}
\end{enumerate}
\end{proposition}
\begin{remark}
\normalfont

    \begin{enumerate}[i)]
    \normalfont
        \item We provide the optimal terminal wealth for a bounded non-concave utility maximization problem and a non-bounded non-concave utility maximization problem, respectively. A non-bounded utility $U(0)=-\infty$ leads to a trivial case mathematically, where the optimal wealth is always above the debt level $D_T$. In this case, the company will not be motivated to promise the debt holders a return that is higher than the risk-free return ($g(t)=r$) since the default risk is eliminated. Economically, this can be used to model a very risk averse attitude of the managerial board which possibly does not allow for default due to reputation costs, etc.
        \item In the non-trivial case ($U(0)>-\infty$), the surplus-driven company distinguishes the good and bad financial states according to the critical value of the state price density $\xi_B$ determined by its risk attitude (parametric utility functions) and the terminal debt. Consequently, the benchmark wealth is either $\wh D_T$, or is zero.\footnote{We remark that $\wh D_T$ is uniquely determined by the utility function and the debt level, and $\wh D_T>D_T$ always holds.} Compared with a non-surplus-driven company that maximizes the expected utility of the total asset, a surplus-driven company invests more in the good financial states to obtain a higher surplus, and ends up with zero terminal wealth in the bad financial states (Figure \ref{poffB}).\footnote{The optimal terminal wealth of a non-surplus-driven company is a direct application of Merton solution \citep{merton1969lifetime}.}
        \item While we consider a simple non-linear payoff function, the benchmark solution and its implication is general enough, and is consistent with the literature. The key intuition is that if the punishment for incurring the zero terminal wealth (or the maximum loss) is not large enough ($U(0)>-\infty$), the company will anyway incur zero terminal wealth in the worst scenarios. If the punishment is too large ($U(0)=-\infty$), the model admits a trivial solution mathematically. However, a company's attitude towards the worst financial scenarios can only be answered by empirical studies. We include both cases for completeness.
    \end{enumerate}
    \end{remark}
		\subsection{A \emph{Concave} utility maximization problem under a risk constraint}
	Since Merton \cite{merton1969lifetime} and Samuelson \cite{samuelson1975lifetime}, optimal asset allocation problems based on expected utility theory have been extensively studied. One important stream of extension considers the optimization problem with an additional risk constraint representing the risk-based regulation. In a complete financial market, explicit solutions are obtained. We find that the optimal solutions for the counterpart concave utility maximization problem under the risk constraints ESP \eqref{def:es}, ESQ \eqref{def:esq}, VaR \eqref{def:var} and AVaR \eqref{def:AVaR} are similar in structure, and these risk constraint mainly affect the tail distribution of the optimal terminal wealth. In order to save space, we do not repeat the explicit solution in each case, but summarize the general structure of the optimal solutions in the following proposition.\footnote{We provide the structure of optimal solutions assuming that the risk constraints are not redundant.}
	
		\begin{proposition}\label{proposition:concave}
	\begin{enumerate}[i)] 
	    \item Let $\lambda_1$ and  $\lambda_2$ be the corresponding Lagrangian multipliers for the budget and one of the risk (ESP,ESQ,VaR) constraints, respectively. The general structure of the optimal solution under one of the (ESQ, ESP, VaR) risk constraints is given by
	\begin{equation}
	 X_T^{*}=(U^{'})^{-1}(\lambda_1\xi)\mathbbm{1}_{\xi<\xi_1}+L\mathbbm{1}_{\xi_1\leq\xi<\xi_2}+f(\lambda_1,\lambda_2,\xi)\mathbbm{1}_{\xi\geq\xi_2},   \end{equation}
	 where $f(\lambda_1,\lambda_2,\xi)$ is a function of the state price density and the Lagrangian multipliers. Its concrete form depends on the exact risk constraint. In addition, $\xi_1$ and $\xi_2$ are two critical values of the state price density, and are also determined by the exact risk constraint.\footnote{One can find the specific form for $f(\lambda_1,\lambda_2,\xi)$ under a VaR or an ESQ constraint in \cite{basak2001value}, under an ESP constraint in \cite{shi2012short}.}
	 \item Under the AVaR constraint, the structure of the optimal solution is given by
	 \begin{equation}
	     X_T^{*}=(U^{'})^{-1}(\lambda_1\xi)\mathbbm{1}_{\xi<\xi_1}+L^{*}\mathbbm{1}_{\xi_1\leq\xi<\xi_2}+f(\lambda_1,\lambda_2,\xi)\mathbbm{1}_{\xi\geq\xi_2},  
	 \end{equation}
	 where $L^{*}$ is a constant but is implicitly determined by the model parameters.\footnote{One can check \cite{wei2021risk} for $L^{*}$ and $f(\lambda_1,\lambda_2,\xi)$ under an AVaR constraint.}
	\end{enumerate}
	\end{proposition}
	\begin{remark}
	\normalfont
	\begin{enumerate}[i)]
    \item The optimal wealth under a risk constraint distinguishes the good, intermediate and bad financial states according to the critical values of the state price density $\xi_1$ and $\xi_2$. The exact values of $\xi_1$ and $\xi_2$ are determined by the initial wealth, the risk attitude of the company, and the concrete risk constraint. Hence, they differ in different solutions. 
    \item The constrained optimal wealth is similar to the unconstrained Merton solution \citep{merton1969lifetime} in the good financial states, is a constant in the intermediate financial states, and is different in the bad financial states under different risk constraints.
    \item By examining the explicit solutions, one can see that the unconstrained optimal wealth in the bad financial states is always lower than the optimal wealth under the ESQ constraint, but is higher than the optimal wealth under the ESP, VaR or AVaR constraint.
      Hence, compared with the other three risk constraints, the ESQ constraint provides a better protection in the concave utility maximization problem.
\end{enumerate}
	\end{remark} 
	In the next section, we will examine how these risk constraints affect the optimal wealth in the non-concave utility maximization. 
		\section{Non-concave utility maximization under a risk constraint}	
		In this section, we provide the explicit solutions for the non-concave utility maximization problems under the risk constraints (ESP,ESQ,VaR,AVaR). Note that if we add a risk constraint to an optimization problem, we always have to check i) whether the optimal solution exists; ii) whether the risk constraint is redundant.  A risk constraint is redundant if the unconstrained solution already satisfies the risk constraint. Given a risk constraint, this occurs if the initial wealth exceeds an upper bound ($x_0\geq\overline{x}_0$). In this case, the unconstrained solution is also the optimal solution under the risk constraint. On the other hand, the optimal solution under a non-redundant risk constraint exists if the initial wealth is sufficiently large ($x_0\geq x_0^{min}$). Note that the exact boundary values of the initial wealth ($x_0^{min}, \overline{x}_0$) depend on the concrete risk constraint. In the following, we provide the optimal solutions under a non-redundant risk constraint assuming that they exist. After obtaining the non-redundant optimal solutions, calculating the boundary values of the initial wealth is straightforward.
		
			Recall that we formulate four optimization problems in Section 2.2.2, i.e., Problem ESQ, Problem ESP, Problem VaR and Problem AVaR. Surprisingly, the structures of the optimal solutions under the four risk constraints are similar. Hence, we provide the general structure of the optimal solutions in Proposition \ref{proposition:nonconcave}. In addition, we give the explicit solutions under the ESP and the AVaR constraints as examples. All other explicit solutions and the corresponding proofs are put in Appendix B.
		\begin{proposition} \label{proposition:nonconcave}
	Let $\lambda_1$ and $\lambda_2$ denote the Lagrangian multipliers of the budget constraint and the risk constraint. Let $\xi_1$, $\xi_2$ and $\xi_3$ be three critical values of the state price density implicitly determined by the risk constraint, the risk attitude of the company, and the terminal debt level. Recall that $\wh D_T$ is the lower bound of the optimal benchmark wealth in the good states (Proposition \ref{proposition:benchmark}).
	\begin{enumerate}[i)]
	    \item The general structure of the optimal solutions for the non-concave utility maximization problems under the (ESQ, ESP or VaR) constraint is given by
	\begin{footnotesize}    
	\begin{equation}
	    X_T^{*}=\begin{cases}
	   & \left((U^{'})^{-1}(\lambda_1\xi)+D_T\right)\mathbbm{1}_{\xi<\xi_1}+L\mathbbm{1}_{\xi_1\leq\xi<\xi_2}+h(\lambda_1,\lambda_2,\xi)\mathbbm{1}_{\xi_2\leq\xi<\xi_3},\quad\text{if}\quad L>\wh D_T;\\
	    &\left((U^{'})^{-1}(\lambda_1\xi)+D_T\right)\mathbbm{1}_{\xi<\xi_1}+L\mathbbm{1}_{\xi_1\leq\xi<\xi_2},\quad\text{if}\quad \xi_1<\xi_2\quad\text{and}\quad L\leq\wh D_T;\\
	    &\left((U^{'})^{-1}(\lambda_1\xi)+D_T\right)\mathbbm{1}_{\xi<\xi_2},\quad\text{if}\quad\xi_2<\xi_1\quad\text{and}\quad L\leq\wh D_T,
	       \end{cases}
	    \end{equation}
	    	\end{footnotesize}
	    where $h(\lambda_1,\lambda_2,\xi)$ is a function of the Lagrangian multipliers and the state price density. Its form depends on the exact risk constraint.\footnote{One can check Proposition \ref{proposition:esp} for $h(\lambda_1,\lambda_2,\xi)$ under the ESP constraint, and Propositions \ref{solu:esq} and \ref{solu var} for $h(\lambda_1,\lambda_2,\xi)$ under the ESQ and VaR constraints.}
	    \item The general structure of the optimal solution under the AVaR constraint is given by
	    \begin{footnotesize}
	    	\begin{equation}
	    X_T^{*}=\begin{cases}
	   & \left((U^{'})^{-1}(\lambda_1\xi)+D_T\right)\mathbbm{1}_{\xi<\xi_1}+L^{*}\mathbbm{1}_{\xi_1\leq\xi<\xi_2}+h^{AVaR}(\lambda_1,\lambda_2,\xi)\mathbbm{1}_{\xi_2\leq\xi<\xi_3},\quad\text{if}\quad L^{*}>\wh D_T;\\
	    &\left((U^{'})^{-1}(\lambda_1\xi)+D_T\right)\mathbbm{1}_{\xi<\xi_1}+L^{*}\mathbbm{1}_{\xi_1\leq\xi<\xi_2},\quad\text{if}\quad \xi_1<\xi_2\quad\text{and}\quad L^{*}\leq\wh D_T;\\
	    &\left((U^{'})^{-1}(\lambda_1\xi)+D_T\right)\mathbbm{1}_{\xi<\xi_2},\quad\text{if}\quad\xi_2<\xi_1\quad\text{and}\quad L^{*}\leq\wh D_T,
	       \end{cases}
	    \end{equation} 
	     \end{footnotesize}
	    where $L^{*}$ is a constant and is implicitly determined by the model parameters, and $h^{AVaR}(\lambda_1,\lambda_2,\xi)$ is a function of the Lagrangian multipliers and the state price density.
	\end{enumerate}
	\end{proposition}
	\begin{remark}\label{remark:nonconcave}
	\normalfont
		\begin{enumerate}[i)]
		 \item We have three different structures of the optimal wealth for a surplus-driven company under the risk constraints. Mathematically, this is because the highly non-concave Lagrangian can be decomposed into two, or three or four local (concave) Lagrangian on disjoint sets depending on the model parameters. Hence, the global maximum of the Lagrangian could exist on any subset. Moreover, if we fix a risk constraint, each structure of the solution exists within a limited range of the initial wealth. Economically, this means that for a fixed risk constraint, the company chooses different strategies according to its initial wealth. The company with a  concave utility function (Proposition \ref{proposition:concave}) has a single structure of the optimal wealth, while a surplus-driven company has multiple structures of the optimal wealth depending on the initial wealth. Hence, the strategy of a surplus-driven company is more sensitive to the change of its initial wealth. 
		  \item Each structure of the optimal solution has a non-negligible set of zero terminal wealth. Compared with the benchmark solution (Proposition \ref{proposition:benchmark}), these risk constraints \emph{do not} eliminate but can reduce this set in the non-concave optimization problems, indicating the success and failure of these risk-based financial regulations on a surplus-driven financial company.
	    
	    \item If the given threshold $L$ or the implicit threshold $L^{*}$ is smaller or equal to $\wh D_T$, the optimal solutions only differ in critical values of state price density $\xi_1$ and $\xi_2$. Thus, if $\xi_1$ and $\xi_2$ are the same in different solutions, these solutions are essentially equivalent. In fact, we can establish a \emph{one-to-one} relationship between risk constraints such that the optimal solutions under the corresponding risk constraints are the \emph{same}. Note that the equivalence among risk constraints in a concave model does not exist due to the various distributions of the optimal wealth in the tail. We discuss the equivalence among risk constraints in the non-concave model in Section 5.
	\end{enumerate}
		\end{remark}
			Now, we provide the complete explicit solution under the (non-redundant) ESP and AVaR constraints as examples. Other explicit solutions can be found in Appendix B.
			
				\paragraph{Case 1: The optimal solution for Problem ESP}

	Recall that Problem ESP is given by
	\begin{equation*}
	\max\limits_{X_T\in\mathcal{X}}\bbe\big[U\big((X_T-D_T)^+\big)\big],\quad\text{s.t.}\quad \bbe[(L-X_T)^+]\leq\epsilon_\bbp,\quad \bbe[X_T\xi]\leq x_0.
	\end{equation*}
		The optimal wealth of Problem ESP is given in the following proposition. The proof can be found in Appendix B.
	\begin{proposition}\label{proposition:esp}	
		\begin{enumerate}[i)]
			\item If $L\leq D_T$, the optimal wealth is:
				\begin{align}
				X_T^{ESP}=&\left((U^{'})^{-1}(\lambda_{1}\xi)+D_T\right)\mathbbm{1}_{\xi<\xi_{2}} \quad\text{if $\xi_{2}\leq \xi_1$,}\label{26}\\
				X_T^{ESP}=&\left((U^{'})^{-1}(\lambda_{1}\xi)+D_T\right)\mathbbm{1}_{\xi<\xi_1}+L\mathbbm{1}_{\xi_1\leq\xi<\xi_2} \quad\text{if $\xi_1<\xi_2$,}\label{esq25}
				\end{align}
					\noindent where $\xi_1=U^{'}(\wt L-(D_T-L))/\lambda_{1}$, $\wt L$ is the tangent point of $U((X_T-D_T+L)^+)$, $\xi_2$ is defined through $\bbe[L\mathbbm{1}_{\xi\geq\xi_2}]=\epsilon_\bbp$, 
				and $\lambda_1$ is obtained by solving $\bbe[\xi X_T^{ESP}]=x_0$.
			\item If $D_T<L\leq\wh D_T$, the optimal wealth is:
					\begin{align}
				X_T^{ESP}=&\left((U^{'})^{-1}(\lambda_{1}\xi)+D_T\right)\mathbbm{1}_{\xi<\xi_{2}}\quad\text{if $\xi_{2}<\xi_1$,}\label{35}\\
				X_T^{ESP}=&\left((U^{'})^{-1}(\lambda_{1}\xi)+D_T\right)\mathbbm{1}_{\xi<\xi_1}+L\mathbbm{1}_{\xi_1\leq\xi<\xi_2}\quad\text{if $\xi_1\leq\xi_2$,}\label{36}
				\end{align}
					where $\xi_1=U^{'}(L-D_T)/\lambda_{1}$ and $\xi_2$ is defined through $\bbe[L\mathbbm{1}_{\xi\geq\xi_2}]=\epsilon_\bbp$, 
				and $\lambda_1$ is obtained by solving $\bbe[\xi X_T^{ESP}]=x_0$.
		\item If $L>\wh D_T$, the optimal wealth is:
		\begin{footnotesize}
		    	\begin{equation}\label{17}
		X_T^{ESP}=\left((U^{'})^{-1}(\lambda_1\xi)+D_T\right)\mathbbm{1}_{\xi<\xi_1}+L\mathbbm{1}_{\xi_1\leq\xi<\xi_2}+\left((U^{'})^{-1}(\lambda_1\xi-\lambda_2)+D_T\right)\mathbbm{1}_{\xi_2\leq\xi<\xi_3},
				\end{equation}
				\end{footnotesize}
					\noindent where $\xi_1=\frac{ U^{'}(L-D_T)}{\lambda_1}$, $\xi_2=\frac{ U^{'}(L-D_T)+\lambda_2}{\lambda_1}$,
			$\xi_3=\frac{U^{'}(\wh D_T-D_T)+\lambda_2}{\lambda_1}$, $\lambda_1$ and $\lambda_2$ are obtained by solving the equations $\bbe[\xi X_T^{ESP}]=x_0$ and $\bbe[(L-X_T^{ESP})^+]=\epsilon_{\bbp}$.
				\end{enumerate}
\end{proposition}
	The optimal wealth for Problem ESP in Proposition \ref{proposition:esp} shows exactly the general structure of the optimal wealth under a risk constraint provided in Proposition \ref{proposition:nonconcave}. Indeed, the specific values of the state price density $\xi_1$, $\xi_2$ and $\xi_3$ are determined by the ESP constraint, the risk attitude (utility function) and the terminal debt level. Especially, if the regulatory threshold $L$ is smaller than $\wh D_T$, $\xi_3$ and $\lambda_2$ are not involved in the explicit solution.
	We obtain similar conclusions for the optimal solutions under the ESQ and the VaR constraints. The case with the AVaR constraint is different. We introduce the optimal solution under the AVaR constraint in the next example.

\paragraph{Case 2: The optimal solution for Problem AVaR}
Recall that Problem AVaR is given by
	\begin{equation*}
	\max\limits_{X_T\in\mathcal{X}}\bbe[U((X_T-D_T)^+)],
	\end{equation*}
	subject to
	\begin{equation*}
	\frac{1}{\alpha}\int_{0}^{\alpha}VaR_{X_T}(\beta)\d\beta\geq L^{AVaR},\quad \bbe[X_T\xi]\leq x_0.
	\end{equation*}
	The maximization problem is formulated based on the terminal wealth $X_T$, while the AVaR constraint is formulated based on the quantile function of the terminal wealth. Hence, the Lagrangian approach cannot be directly applied. We deal with the problem in two steps: i) We transfer the AVaR constraint to a set of equivalent ESP constraints; ii) We show that the optimal solution under the AVaR constraint is equivalent to the optimal solution under one of the equivalent ESP constraint. The following lemma establishes the connection between the AVaR constraint and the ESP constraint.
\begin{lemma}\label{lemma:AVaR}
    \begin{enumerate}[i)]
        \item If the threshold $L$ in an ESP constraint is exactly the $\alpha$-quantile of the terminal wealth, then the AVaR constraint is equivalent to the following ESP constraint:
        \begin{equation}
            \frac{1}{\alpha}\int_0^{\alpha}VaR_{X_T}(\beta)\d \beta\geq L^{AVaR}\iff\bbe[(L-X_T)^+]\leq\epsilon_\bbp,
        \end{equation}
        where $\epsilon_\bbp=\alpha(L-L^{AVaR}).$
         \item Consider the function $\mathcal{G}(L): L\rightarrow\bbe[(L-X_T)^+]-\alpha(L-L^{AVaR})$. 
         The function $\mathcal{G}(L)$ attains its minimum if $L$ is the $\alpha$-quantile of $X_T$;
         \item A terminal wealth that satisfies a given risk constraint is called a feasible solution. The set of the feasible solutions for an AVaR constraint is equivalent to the set of the feasible solutions for the ESP constraints (regardless of $L$ being the $\alpha$-quantile of $X_T$), i.e.,
         \begin{footnotesize}
                  \begin{equation}
             \bigg\{X_T\vline\quad \frac{1}{\alpha}\int_{0}^{\alpha}VaR_{X_T}(\beta)\d\beta\geq L^{AVaR}\bigg\}=\bigcup_{L\geq L^{AVaR}}\bigg\{X_T\vline\quad\bbe[(L-X_T)^+]\leq\alpha\big(L-L^{AVaR}\big)\bigg\}.
         \end{equation}
          \end{footnotesize}
    \end{enumerate}
   \end{lemma}
   Argument $i)$ in Lemma \ref{lemma:AVaR} follows from Proposition 4.51 in \cite{follmer2016stochastic}. Arguments $ii),iii)$ follow from Lemma 2.7 and Theorem 2.9 in \cite{gandy2005portfolio}. Lemma \ref{lemma:AVaR} reveals the connection between the quantile-based AVaR constraint and the terminal wealth-based ESP constraint. We show in the following proposition that the optimal solution under the AVaR constraint is the same as the optimal solution under the equivalent ESP constraint. The proof of Proposition \ref{proposition:AVaR} can be found in Appendix B. 
 
  \begin{proposition}\label{proposition:AVaR}
   Suppose the AVaR constraint is given by: $\frac{1}{\alpha}\int_{0}^{\alpha}VaR_{X_T}(\beta)\d\beta\geq L^{AVaR}$. Let $X_T^{ESP,L}$ denote the optimal solution for the following Problem $ESP$,
    		\begin{equation*}\label{pro:esq}
	\left[	\max\limits_{X_T\in\mathcal{X}}\bbe[U((X_T-D_T)^+)],\quad\text{s.t.}\quad\bbe[(L-X_T)^+]\leq\alpha(L-L^{AVaR}),\quad \bbe[X_T\xi]\leq x_0\right].
		\end{equation*}
For a given initial wealth $x_0>x_0^{min}$, there exists $L^{*}>L^{AVaR}$ such that the optimal solution for Problem AVaR coincides with the optimal solution for Problem ESP with the constraint: $\bbe[(L^{*}-X_T)^+]\leq \alpha(L^{*}-L^{AVaR})$, i.e.,
$$X_T^{AVaR}=X_T^{ESP,L^{*}}.$$
		  \end{proposition}
		   Lemma \ref{lemma:AVaR} and Proposition \ref{proposition:AVaR} together reveal the connection between the optimal solutions under the AVaR constraint and the ESP constraint. We remark that due to the sophisticated relationship between the various structure of the optimal solution and the initial wealth, a more explicit form of $L^{*}$ cannot be computed within the current model. However, it is sufficient to conclude that the optimal solution for Problem AVaR shows an identical structure as the optimal solution for Problem ESP. 
		   In order to better understand the economical implication of the optimal wealth, we provide a numerical illustration of the optimal solution for Problem ESP in a simple Black-Scholes market.
		   
		   \paragraph{Numerical illustration}  Consider a simple Black-Scholes financial market equipped with one risky asset and one risk-free asset. The basic parameters about the financial market are given in Table \ref{table:basic}. In addition, assume that the surplus-driven company has a power utility function $U(x)=\frac{x^{1-\gamma}}{1-\gamma}$ with a constant relative risk aversion $\gamma=0.5$. Thus, the utility is bounded from below ($U(0)>-\infty$). 
	\begin{table}[ht!]
		\begin{center}
			\begin{tabular}{ccccc}
				\hline
				$\mu$&$r$&$\sigma$&$T$&$\gamma$\\
				0.08&0.03&0.2&1&0.5\\
				\hline
			\end{tabular}
			\caption{In a Black-Scholes market, the risky asset is driven by a Brownian motion. The mean and volatility of the risky asset are denoted by $\mu$ and $\sigma$, respectively. The risk-free rate is constant and is denoted by $r$. We consider the time horizon $T=1$.}	\label{table:basic}
		\end{center}
	\end{table}
	
\begin{figure}[ht!]
	\begin{subfigure}{0.48\textwidth}
		\centering
		\includegraphics[scale=0.3]{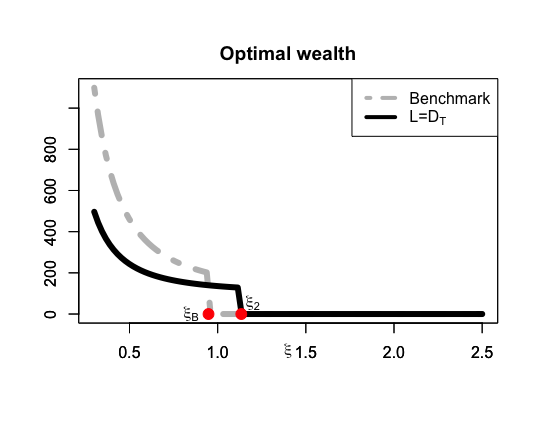}%
		\vspace{-10mm}
		\captionsetup{justification=centering}
		\caption{Case 1}
	\end{subfigure}
	\begin{subfigure}{0.48\textwidth}
		\centering
		\includegraphics[scale=0.3]{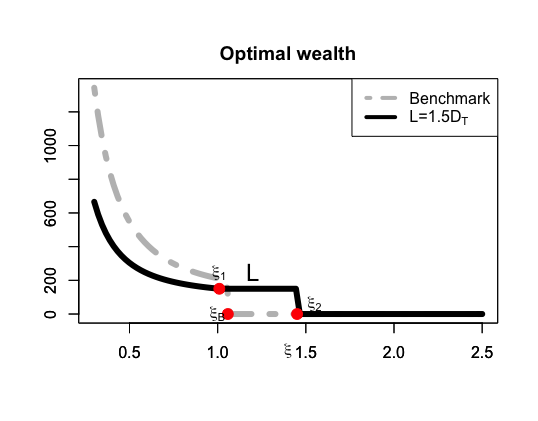}%
		\vspace{-10mm}
		\captionsetup{justification=centering}
				\caption{Case 2}
	\end{subfigure}
		\begin{subfigure}{0.48\textwidth}
		\centering
		\includegraphics[scale=0.3]{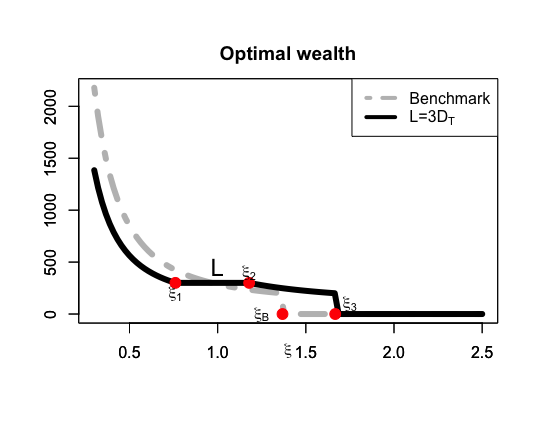}%
		\vspace{-10mm}
		\captionsetup{justification=centering}
				\caption{Case 3}
	\end{subfigure}
	\captionof{figure}{The three different structures of the optimal wealth for the non-concave utility maximization under an ESP constraint. The parameters in each case correspond to the three cases in Table \ref{table:structure}. }\label{fig:ESQ}
\end{figure}

	In addition, to illustrate different structures of the optimal wealth, we consider different combinations of model parameters (e.g. initial wealth and ESP constraints) in Table \ref{table:structure}. Note that the debt level is assumed to be 100 in each case, i.e., $D_T=100$. Consequently, $\wh D_T$, the lower bound of the benchmark wealth in the good financial states (Proposition \ref{proposition:benchmark}) is also fixed in each case.\footnote{Mathematically, $\wh D_T$ is the unique tangent point of $U((X_T-D_T)^+)$, and is determined by the utility function and the debt level (Lemma \ref{tang}).} 
	\begin{table}[ht!]
		\begin{center}
			\begin{tabular}{cccccc}
				\hline
				&$D_T$&$\wh D_T$&$x_0$&$L$&$\epsilon_\bbp$\\
				\hline
				Case 1:&100&200&100&100&30\\
				Case 2:&100&200&150&150&10\\
				Case 3:&100&200&292&300&7\\
				\hline
			\end{tabular}
			\caption{This table contains three potential combinations of the model parameters, e.g. the initial wealth and the different ESP constraints to illustrate the three structures of the optimal wealth (Proposition \ref{proposition:esp}).}	\label{table:structure}
		\end{center}
	\end{table}
		We plot the optimal wealth in the three different cases in Figure \ref{fig:ESQ}. 
		
		From Figure \ref{fig:ESQ} we observe the following: i) the structure of the optimal wealth depends on the model parameters as explained in Remark \ref{remark:nonconcave}; ii) the optimal wealth is zero in the worst financial states in all cases. However, the set of the zero wealth is always smaller than the corresponding benchmark (unconstrained) wealth.
		
		We obtain a similar conclusion for the optimal wealth under the other risk constraints as well (Propositions \ref{solu:esq} and \ref{solu var}). Hence, we conclude that these risk constraints (ESQ, ESP, VaR, AVaR) have a similar regulatory effect, that is to \emph{reduce the set of zero wealth} of a surplus-driven company.

In fact, if the threshold $L$ is below the minimum benchmark wealth in the good financial states, i.e., if $L\leq\wh D_T$, the optimal wealth under different risk constraints (ESQ, ESP, VaR) only differ in the critical values of the state price density $\xi_1$ and $\xi_2$. This implies that these optimal wealth are the \emph{same} if $\xi_1$ and $\xi_2$ are the \emph{same} in different solutions. Moreover, if the implicit threshold $L^{*}$ is below $\wh D_T$, the equivalence result also holds for the optimal wealth under the AVaR constraint. Note that the equivalence holds in a \emph{complete} financial market and does not depend on a specific utility function or a Black-Scholes market assumption. In the next section, we investigate the equivalence among the risk constraints.
\section{Equivalence among risk constraints in non-concave optimization}
In this section, we discuss a by-product of obtaining the optimal solutions for the non-concave utility maximization problems under the risk constraints. We know from Proposition \ref{proposition:nonconcave} that, if the given threshold $L$ or the implicit threshold $L^{*}$ is below $\wh D_T$, the optimal solutions only differ in the critical values of the state price density. Hence, the optimal wealth under these risk constraints are the \emph{same} if $\xi_1$ and $\xi_2$ are the \emph{same}. We summarize this finding in Corollary \ref{col:equii}.

\begin{corollary}\label{col:equii}	\normalfont
 Let $X_T^{ESQ}$, $X_T^{ESP}$, $X_T^{VaR}$ and $X_T^{AVaR}$ denote the optimal wealth for the non-concave uitility maximization problem under the ESQ, ESP, VaR and AVaR constraints, respectively. 
	\begin{enumerate}[i)]
		\item For each given VaR constraint with a significance level $\alpha$ \eqref{def:var}, there exists one ESQ constraint with $\epsilon_\bbq(\alpha)$ and one ESP constraint with $\epsilon_\bbp (\alpha)$ such that $X_T^{VaR}=X_T^{ESQ}=X_T^{ESP}$ if $L\leq\wh D_T$. Note that with the same regulatory threshold $L$, the equivalent (ESP and ESQ) constraint is unique. 
		\item For each given AVaR constraint \eqref{def:AVaR}, there exists one VaR constraint $\bbp(X_T\geq L^{*})\geq 1-\alpha(L^{*})$, one ESQ constraint $\bbe[\xi(L^{*}-X_T)^+]\leq\epsilon_\bbq(L^{*})$ and one ESP constraint $\bbe[(L^{*}-X_T)^+]\leq\epsilon_\bbp(L^{*})$ such that $X_T^{AVaR}=X_T^{ESQ}=X_T^{ESP}=X_T^{VaR}$. The equivalent (ESP, ESQ and VaR) constraint is unique.
	\end{enumerate}
\end{corollary}
\proof
By examining the explicit optimal wealth under different risk constraints in Section 4 and Appendix B, one can obtain this conclusion via direct calculation.\endproof

\noindent Next, we use an numerical example in a Black-Scholes market to illustrate the equivalence. 

\paragraph{Example} Consider the same simple Black-Scholes market with parameters in Table \ref{table:basic}. 
The terminal debt for the surplus-driven company (Section 4) is fixed, i.e., $D_T=100$. Then, the benchmark (unconstrained) wealth is either higher than $\wh D_T=200$ or is zero.

Consider the AVaR constraint 
	$$\frac{1}{\alpha}\int_{0}^{\alpha}VaR_{X_T}(\beta)\d\beta\geq 0.9D_T,\quad\alpha=0.01.$$ 
Thus, the AVaR-based regulation requires that the average quantile of the company's wealth in the worst 1\% scenarios is at least 90\% of its terminal debt.

Recall that (Proposition \ref{proposition:nonconcave}) the optimal wealth under the AVaR constraint is given by
\begin{equation*}
	    X_T^{*}=\begin{cases}
	    &\left((U^{'})^{-1}(\lambda_1\xi)+D_T\right)\mathbbm{1}_{\xi<\xi_1}+L^{*}\mathbbm{1}_{\xi_1\leq\xi<\xi_2},\quad\text{if}\quad \xi_1<\xi_2,\\
	    &\left((U^{'})^{-1}(\lambda_1\xi)+D_T\right)\mathbbm{1}_{\xi<\xi_2},\quad\text{if}\quad\xi_2<\xi_1.
	       \end{cases}
	    \end{equation*} 

The threshold $L^{*}$ is implicitly determined by the parameters. We construct two cases with different $L^{*}$ illustrating different structures of optimal wealth.
\paragraph{Case 1: $L^{*}=L=D_T$.}
The optimal wealth under the AVaR constraint is given by
$$X_T^{*}=\left((U^{'})^{-1}(\lambda_1\xi)+D_T\right)\mathbbm{1}_{\xi<\xi_2}.$$
In the given BS market, we calculate the corresponding VaR, ESQ and ESP constraints such that the optimal wealth has the same structure and $\xi_2$ is the same for each optimal wealth. To be more precise, the equivalent VaR constraint is
$$\bbp(X_T\geq L^{*})\geq99.9\%.$$
Hence, the equivalent VaR constraint requires that the company's wealth is above the threshold $L^{*}=D_T$ with probability 99.9\%. The equivalent ESP constraint is 
$$\bbe[(L^{*}-X_T)^+]\leq 0.1\%D_T,$$
which requires that the expected shortfall measured by $L^{*}$ is no more than 0.1\% of its terminal debt. The equivalent ESQ constraint is 
$$\bbe[\xi(L^{*}-X_T)^+]\leq 0.22\% D_T,$$
which requires that the expected discounted shortfall measured by $L^{*}$ is no more than 0.22\% of its terminal debt level.

\paragraph{Case 2: $L^{*}=L=1.2D_T$.}
The optimal wealth under the AVaR constraint is
$$X_T^{*}=\left((U^{'})^{-1}(\lambda_1\xi)+D_T\right)\mathbbm{1}_{\xi<\xi_1}+L^{*}\mathbbm{1}_{\xi_1\leq\xi<\xi_2}.$$
Again, we calculate the equivalent VaR, ESP, and ESQ constraint such that the optimal wealth has the same structure and $\xi_1$ and $\xi_2$ are the same in different solutions. The equivalent VaR constraint is
$$\bbp(X_T\geq L^{*})\geq99.7\%,$$
which requires that the company's wealth is above the threshold $L^{*}=1.2D_T$ with probability 99.7\%. The equivalent ESP constraint is 
$$\bbe[(L^{*}-X_T)^+]\leq 0.3\%D_T,$$
which requires that the expected shortfall measured by $L^{*}$ is no more than 0.3\% of its terminal debt. The equivalent ESQ constraint is 
$$\bbe[\xi(L^{*}-X_T)^+]\leq 0.73\% D_T,$$
which requires that the expected discounted shortfall measured by $L^{*}$ is no more than 0.73\% of its terminal debt level.

In this section, we have shown an interesting by-product that the optimal wealth under different risk constraints can be made equivalent, if $L$ or $L^{*}$ is smaller than $\wh D_T$. Recall that the benchmark (unconstrained) wealth is either higher than $\wh D_T$ or is zero. Considering the fact that $\wh D_T>D_T$ always holds, the case $L\leq\wh D_T$ is more relevant in practice, since the regulatory threshold is always set close to the debt level. This means that in most cases all the four risk constraints have the same regulatory effect on a surplus-driven financial company. There is no equivalence result in the counterpart concave maximization problems (Proposition \ref{proposition:concave}), since the optimal wealth behaves very differently in the bad financial states under different risk constraints. However, if the company is surplus-driven, these risk constraints are as good as or as bad as each other.

\section{Conclusion}
This paper investigates the effect of risk constraints in non-concave optimization. We solve the optimization problems with a non-concave utility function under four prevalent risk constraints ESP, ESQ, VaR and AVaR, by the Lagrangian approach. Due to the non concavity and the existence of a risk constraint, we obtain several structures of the optimal solutions depending on the model parameters. Although these risk constraints are widely discussed in academia and in the industry, a comparative analysis of these risk constraints on the non-concave utility maximization problem based on the closed form solutions is missing in the literature. Our results fill this gap. 

Moreover, we find that the four risk constraints, though defined differently, all lead to the same optimal wealth in the non-concave optimization problems. They equivalently and effectively reduce the set of zero terminal wealth for a surplus-driven company. Hence, we conclude that considering the potential surplus-driven characteristic of a financial company, the four risk constraints are as good as or as bad as each other. 

Our results imply that it is not sufficient to rely on the expectation or quantile-based risk measure to restrict the gambling behavior of a surplus-driven financial company. A natural research question is what could be a better substitute risk measure. We leave this topic for future research.
\bibliography{22nd}
	\appendix
	\section{Using Lagrangian to solve the optimization problem with a non-concave utility function under risk constraints}
In this section, we briefly summarize the procedure for solving the optimization problems by the Lagrangian approach.

\paragraph{Step 1: Constructing the Lagrangian.} Consider the non-concave objective function $U((X-D_T)^+)$. Denoting by $\rho(X_T)$ one of the risk constraints, we write the budget constraint and the risk constraint as inequalities: $\bbe[\xi X_T]\leq x_0$ and $\rho(X_T)\leq \epsilon$. For instance, the VaR constraint can be written as $\bbp(X_T<L)\leq\alpha$. Letting $\lambda_1$ and $\lambda_2$ be two non-negative Lagrangian multipliers, the constructed Lagrangian is given by:
\begin{equation}\label{lagrangian}
    \phi(X_T):=U((X_T-D_T)^+)-\lambda_1(\xi X_T-x_0)-\lambda_2(\rho(X_T)-\epsilon).
\end{equation}
Note that if the constraints are satisfied by strict inequalities, the Lagrangian multipliers are zero and the Lagrangian multipliers are positive if the constraints are satisfied by equalities. Hence, the optimal solution for the Lagrangian is the same as the optimal solution for the original problem. Note that strict convexity of the utility function is not required for this construction, and thus the Lagrangian can still be analysed if the utility is not globally concave.

\paragraph{Step 2: Determining the global maximum of the Lagrangian.} Due to the non-concave function and the risk constraints, the Lagrangian (\ref{lagrangian}) is highly non-concave. But if we observe closely, the Lagrangian has sufficient convexity such that the problem can be solved. We summarize the method in the following Lemma.
\begin{lemma}\label{shit}
    Suppose there exists a partition of $\Omega$, i.e., $\bigcup_{i\in I}A_{i}=\Omega$ and $\forall i,j\in I, A_i\cap A_j=\emptyset$ such that the local Lagrangian $\phi_i(X_T), X_T\in A_i$ is strictly concave, convex or affine. Letting $X_T^{i}$ denote the argmax of the local Lagrangian obtained by the Kuhn Tucker condition, i.e., $\phi_i(X_T^i)\geq \phi_i(X_T^i), \forall X_T\in A_i$, then the optimal solution is given by
    \begin{equation}
        X_T^{*}=X_T^{i}(\omega), \quad\mbox{if}\quad\omega\in\{\phi_i(X_T^i)>\phi_j(X_T^j)\forall j\in I\}\bigcap A_i.
    \end{equation}
\end{lemma}
Lemma \ref{shit} follows directly from the generalised convex programs (\cite{rockafellar2015convex}). We only need to decompose the non-concave Lagrangian into local concave or affine Lagrangian on disjoint sets, and the global maximum of the Lagrangian is obtained by comparing the local maximums of the Lagrangian. 

After we decompose the Lagrangian, we find that determining the global maximum reduces to solving a sequence of \emph{conjugate functions.} The so-called conjugate functions essentially arise from comparing the local maximums of Lagrangian. We find that the optimal solutions are always zero on the set where the conjugate functions are negative. Intuitively, it means that one should not invest if the obtained utility is smaller than the cost for satisfying the constraints.

	\begin{definition}
		\normalfont For a fixed and positive $d$, consider a concave utility function $U(x-d)$ defined on $(d,\infty)$ which satisfies the Inada and AE condition.  For a given positive $\lambda$, we define the conjugate function $c(y):=\sup\limits_{x>d}\{U(x-d)-x\lambda y\}=U(x(y)-d)-U^{'}(x(y)-d)x(y)=U(I(\lambda y))-\lambda yI(\lambda y)-\lambda yd$, where $I$ is the inverse function of the first derivative of $U$ and $y>0$.
	\end{definition}
	
\begin{figure}[ht]
		\centering
		\includegraphics[width=0.5\textwidth]{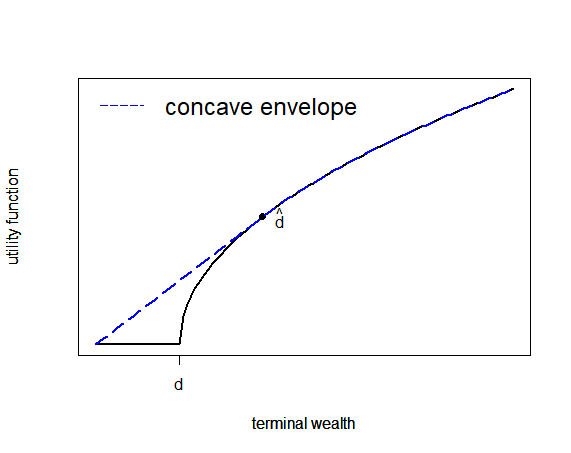}%
		\caption{The non-concave utility function and the corresponding concave envelope.}
		\label{cone}
	\end{figure}

		\begin{lemma} \label{tang}
		\normalfont 
		\begin{enumerate}[i)]
			\item $c(y)$ is decreasing in $y$.
			
			\item  For each positive and fixed $d$, there exists a $\hat{d}$ such that $U(\hat{d}-d)/\hat{d}=U^{'}(\hat{d}-d)$. We call $\hat d$ the \emph{tangent point} of function $U(x-d)$. We plot the tangent point $\wh{d}$ and the concave envelope of function $U(x-d)$ in Figure \ref{cone}. In addition, we have that $y(\hat{d}):=U^{'}(\hat{d}-d)/\lambda$ and $c(y(\hat{d}))=0$. In other words, $y(\hat{d})$ is the unique \emph{zero root} of the conjugate function $c(y)$.
			
			\item If $y>y(\hat{d})$, we have that $U(x(y)-d)/x(y)-U^{'}(x(y)-d)<0$ where $x(y)\equiv I(\lambda y)+d$.
			
			\item For positive $\lambda,\lambda_2,l$ and $d$, we define a conjugate function $c^{*}(y):=\sup\limits_{x>d}\{U(x-d)-x\lambda y+\lambda y\frac{\lambda_2 l}{\lambda}\}=U(I(\lambda y))-\lambda yI(\lambda y)-\lambda yd+\lambda y\frac{\lambda_2 l}{\lambda}$. Then $c^{*}(y)$ has a zero root if and only if $\frac{\lambda_2 l}{\lambda}-d=:s<0$.
				\item For positive $\lambda,\lambda_3,l$ and $d$, we define a conjugate function $c^{**}(y):=\sup\limits_{x>d}\{U(x-d)-x\lambda y+\lambda_3 l\}=c(y)+\lambda_3 l$. Then $c^{**}(y)$ has a zero root $y^{**}$, and $y^{**}>y(\hat{d})$.
			\item We define a conjugate function $c^{***}(y):=\sup\limits_{x>d}\{U(x+d)-x\lambda y\}=U(x(y)+d)-U^{'}(x(y)+d)x(y)=U(I(\lambda y))-\lambda yI(\lambda y)+\lambda yd$, where $x(y)=I(\lambda y)-d$, decreases on $(0,U^{'}(d)/\lambda)$, increases on $(U^{'}(d)/\lambda,\infty)$ and hence has the minimum value $U(d)$.
		\end{enumerate}	
	\end{lemma}
		\proof
		The conjugate function describes the \textit{least upper hyperplane} of the static Lagrangian $U(x-d)-x\lambda y$. For a fixed and positive $y$, we have that $c(y):=U(x(y)-d)-U^{'}(x(y)-d)x(y)=U(I(\lambda y))-\lambda yI(\lambda y)-\lambda yd$ since $U^{'}(x(y)-d)\equiv\lambda y$ holds for the conjugate function, denoting by $I$ the inverse function of the first derivative of $U$. The derivative of $c(y)=U(x(y)-d)-\lambda yx(y)$ with respect to $y$ is $-\lambda x(y)$. As $x(y)\equiv I(\lambda y)+d>0$, statement i) in Lemma \ref{tang} follows. 
		
		Note that the utility function satisfies the following Inada and asymptotic elasticity (AE) condition:
		$$\textbf{Inada}:\lim_{x\to 0}U^{'}(x)=\infty;\quad\lim_{x\to\infty}U^{'}(x)=0.$$
		$$\textbf{AE}:\lim\limits_{x\to\infty}\sup\frac{xU^{'}(x)}{U(x)}<1.$$
		Let $I(\lambda y)=z$ and rewrite $c(y)$ as $U(z)-U^{'}(z)z-U^{'}(z)d$. From the Inada and the AE condition we have that
		\begin{align}\label{y0}
		\lim_{z \to \infty}U(z)-U^{'}(z)z-U^{'}(z)d&=\lim_{z \to \infty}U(z)\left(1-\frac{U^{'}(z)z}{U(z)}\right)-U^{'}(z)d>0,
		\end{align}
		\begin{align}\label{yinf}
		\lim_{z \to 0}U(z)-U^{'}(z)z-U^{'}(z)d&=\lim_{z \to 0}U(z)\left(1-\frac{U^{'}(z)z}{U(z)}\right)-U^{'}(z)d<0.
		\end{align}
		In addition, $c(y)$ is a continuous function. Hence by the intermediate value theorem, for each positive $d$, $c(y)$ has a zero root. Thus, statement ii) in Lemma \ref{tang} follows.
		
			We see that $c^{**}(y)=c(y)+\lambda_3 l=U(z)-U^{'}(z)z-U^{'}(z)d+\lambda_3 l$. With a similar argument, $c^{**}(y)$ has a unique zero root denoted by $y^{**}$. Since $c(y)$ and $c^{**}(y)$ are decreasing functions in $y$ (by statement i)), we obtain that $y^{**}>y(d)$. Thus, statement v) is shown.
			
				Statement iii) is a direct application of statement i) and ii). Statement vi) is also a consequence of statement i). Now we only need to show statement iv)
		
		We have $\partial c^{*}(y)/\partial y=\lambda(\frac{\lambda_2 l}{\lambda}-x)$ and $x=I(\lambda y)+d$ holds. We can see that if $s:=\frac{\lambda_2 l}{\lambda}-d<0$, $\lambda(\frac{\lambda_2 l}{\lambda}-x)=\lambda(s-I(\lambda y))<0$. Therefore, $h$ is a decreasing function in $y$. We can write the conjugate function as $U(I(\lambda y))-\lambda y(I(\lambda y)-s)$. It follows from statement ii) in Lemma \ref{tang} that $c^{*}(y)$ has a unique zero root. \endproof
		
			Lemma \ref{tang} covers all the conjugate functions we encountered in solving the non-concave optimization problems under risk constraints. Hence, with Lemma \ref{tang}, it is enough to observe optimal solutions for Problems ESQ, ESP, VaR and AVaR. The proofs are in Appendix B. 
				\paragraph{Step 3: Showing the existence of the Lagrangian multipliers.} In the previous steps, we construct the Lagrangian and solve the problem assuming that the Lagrangian multipliers exist. The last step is to show that the Lagrangian multipliers indeed exist such that the Lagrangian method works. We postpone this proof to Appendix C.
				
						\section{The optimal solutions for Problems ESQ and VaR}
		In this section, we derive the optimal solutions for Problems ESQ and VaR. The proof consists of three steps: 1) We use two Lagrangian multipliers to capture the budget constraint and the risk constraint and establish the static Lagrangian. 2) Then, we find out the global maximum of the Lagrangian by discussing the zero roots of the corresponding conjugate functions. 3) We show that the maximizer (where the global maximum is achieved) is indeed the optimal solution to the problem. Because the proofs for different optimal solutions are similar, we only provide the most difficult case, the ESQ problem, in detail. 
		\subsection{The optimal solution for Problem ESQ}\label{sec:proofesq}
			The optimal wealth of Problem ESQ is given in the following proposition. 
				\begin{proposition}\label{solu:esq}	
		\normalfont 
		\begin{enumerate}[i)]
			\item If $L\leq D_T$, the optimal wealth is:
				\begin{align}
				X_T^{ESQ}=&\left((U^{'})^{-1}(\lambda_{1}\xi)+D_T\right)\mathbbm{1}_{\xi<\xi_{2}} \quad\text{if $\xi_{2}\leq \xi_1$,}\label{26}\\
				X^{ESQ}=&\left((U^{'})^{-1}(\lambda_{1}\xi)+D_T\right)\mathbbm{1}_{\xi<\xi_1}+L\mathbbm{1}_{\xi_1\leq\xi<\xi_2} \quad\text{if $\xi_1<\xi_2$,}\label{esq25}
				\end{align}
					\noindent where $\xi_1=U^{'}(\wt L-(D_T-L))/\lambda_{1}$, $\wt L$ is the tangent point of $U((X_T-D_T+L)^+)$, $\xi_2$ is defined through $\bbe[\xi L\mathbbm{1}_{\xi\geq\xi_2}]=\epsilon_\bbq$, 
				and $\lambda_1$ is obtained by solving $\bbe[\xi X_T^{ESQ}]=x_0$.
						\item If $D_T<L\leq\wh D_T$, the optimal wealth is:
							\begin{align}
				X_T^{ESQ}=&\left((U^{'})^{-1}(\lambda_{1}\xi)+D_T\right)\mathbbm{1}_{\xi<\xi_{2}}\quad\text{if $\xi_{2}<\xi_1$,}\label{35}\\
				X_T^{ESQ}=&\left((U^{'})^{-1}(\lambda_{1}\xi)+D_T\right)\mathbbm{1}_{\xi<\xi_1}+L\mathbbm{1}_{\xi_1\leq\xi<\xi_2}\quad\text{if $\xi_1\leq\xi_2$,}\label{36}
				\end{align}
					where $\xi_1=U^{'}(L-D_T)/\lambda_{1}$ and $\xi_2$ is defined through $\bbe[\xi L\mathbbm{1}_{\xi\geq\xi_2}]=\epsilon_\bbq$, 
				and $\lambda_1$ is obtained by solving $\bbe[\xi X_T^{ESQ}]=x_0$.
		\item If $L>\wh D_T$, the optimal wealth is:
		\begin{footnotesize}
		         	\begin{equation}\label{17}
		X_T^{ESQ}=\left((U^{'})^{-1}(\lambda_1\xi)+D_T\right)\mathbbm{1}_{\xi<\xi_1}+L\mathbbm{1}_{\xi_1\leq\xi<\xi_2}+\left((U^{'})^{-1}(\xi(\lambda_1-\lambda_2))+D_T\right)\mathbbm{1}_{\xi_2\leq\xi<\xi_3},
				\end{equation}
				\end{footnotesize}
					\noindent where $\xi_1=\frac{ U^{'}(L-D_T)}{\lambda_1}$, $\xi_2=\frac{ U^{'}(L-D_T)}{\lambda_1-\lambda_2}$,
			$\xi_3=\frac{U^{'}(\wh D_T-D_T)}{\lambda_1-\lambda_2}$, $\lambda_1$ and $\lambda_2$ are obtained by solving the equations $\bbe[\xi X_T^{ESQ}]=x_0$ and $\bbe[(L-X_T^{ESQ})^+]=\epsilon_{\bbq}$. 
				\end{enumerate}
	\end{proposition}
		\proof 
			\paragraph{Case 1: $L\leq D_T$.}
			The Lagrangian is given by
		\begin{align*}
		\phi(X_T)&= U((X_T-D_T)^+)-\lambda_{1}\xi_TX_T-\lambda_2\xi_T(L-X_T)\mathbbm{1}_{X_T<L}\\
		&=\begin{cases}
		U(X_T-D_T)-\lambda_{1}\xi_TX_T&\quad\text{if}\quad X_T>D_T,\nonumber\\
		-\lambda_{1}\xi_TX_T-\lambda_2\xi_T(L-X_T)\mathbbm{1}_{X_T<L}&\quad\text{if}\quad 0 \leq X_T \leq D_T.
		\end{cases}
		\end{align*}
			Letting $I()$ denote the inverse function of the first derivative of the utility function, i.e., $I(\cdot)=(U^{'})^{-1}(\cdot)$, the first part of the Lagrangian attains its maximum at $I(\lambda_{1}\xi_T)+D_T$. The second part is a negative affine function and obtains its maximum values at 0 or at the jump point $L$ since $\lambda_2<\lambda_{1}$. Further, we have $\phi(0)=-\lambda_2\xi_TL >\phi(L)=-\lambda_{1}\xi_TL$. Therefore, 0 is the local maximizer (argmax of the maximum) of the second part. Moreover, we can see that
		\begin{align}\label{esq1}
		\phi(I(\lambda_{1}\xi_T)+D_T)-\phi(0)&=U(I(\lambda_{1}\xi_T))-\lambda_{1}\xi_TI(\lambda_{1}\xi_T)-\lambda_{1}\xi_TD_T+\lambda_2\xi_TL\nonumber\\
		&=\max\limits_{X_T} \bigg\{U(X_T-D_T)-X_T\lambda_{1}\xi_T+\lambda_{1}\xi_T\frac{\lambda_2}{\lambda_{1}}L\bigg\}.
		\end{align}
		Lemma \ref{tang} gives the necessary and sufficient condition for the existence of the zero root of (\ref{esq1}). In addition, if $L\leq D_T$, we have that $\lambda_{1}\big(\frac{\lambda_2}{\lambda_{1}}L-D_T\big)<\lambda_{1}(L-D_T)\leq0$. Therefore, if $\lambda_2<\lambda_{1}$, there is a zero root of (\ref{esq1}). 
			In particular, if $\lambda_2$ tends to 0, (\ref{esq1}) uniformly converges to $\max\limits_{X_T} \{U(X_T-D_T)-X_T\lambda_{1}\xi_T\}$ (of which the zero root is $U^{'}(\wh D_T-D_T)/\lambda_1$) on each compact set, and thus the zero root of (\ref{esq1}) converges to $U^{'}(\wh D_T-D_T)/\lambda_1$. If $\lambda_2$ tends to $\lambda_{1}$,  (\ref{esq1}) uniformly converges to $\max\limits_{X_T} \{U(X_T-D_T)-X_T\lambda_{1}\xi_T+\lambda_{1}\xi_TL\}$( of which the zero root is $U^{'}(\wt L-(D_T-L))/\lambda_1$), and the zero root of (\ref{esq1})
			converges to $U^{'}(\wt L-(D_T-L))/\lambda_1$. 
				In the case $\xi_{2}>\xi_1$, we set $\lambda_2=\lambda_{1}$. Then, the affine part of the Lagrangian is constant on $X_T<L$. This means that both $L$ and 0 can be the local maximizer in the affine part. Moreover, 
		\begin{align}\label{27}
		\phi(I(\lambda_{1}\xi_T)+D_T)-\phi(L)&=U(I(\lambda_{1}\xi_T))-\lambda_{1}\xi_TI(\lambda_{1}\xi_T)-\lambda_{1}\xi_TD_T+\lambda_{1}\xi_TL\nonumber\\
		&=\max\limits_{X_T} \{U(X_T-D_T)-X_T\lambda_{1}\xi_T+\lambda_{1}\xi_TL\}.
		\end{align}
			We know from Lemma \ref{tang} that (\ref{27}) has a unique zero root which is denoted by $U^{'}(\wt L-(D_T-L))/\lambda_1$. Hence, (\ref{esq25}) is the argmax of the Lagrangian.
		
		\paragraph{Case 2: $L>D_T$}
			The static Lagrangian in this case is
		\begin{align*}
		\phi(X_T)&= U((X_T-D_T)^+)-\lambda_{1}\xi_TX_T-\lambda_2\xi_T(L-X_T)\mathbbm{1}_{X_T<L}\\
		&=\begin{cases}
		U(X_T-D_T)-\lambda_{1}\xi_TX_T&\quad\text{if}\quad L<X_T,\nonumber\\
		U(X_T-D_T)-\lambda_{1}\xi_TX_T-\lambda_2\xi_T(L-X_T)&\quad\text{if}\quad D_T<X_T\leq L,\nonumber\\
		-\lambda_{1}\xi_TX_T-\lambda_2\xi_T(L-X_T)&\quad\text{if}\quad 0\leq X_T \leq D_T.
		\end{cases}
		\end{align*}
		
			The Lagrangian has four possible local maximizers: $I(\lambda_{1}\xi_T)+D_T$ on the first part, $I((\lambda_{1}-\lambda_2)\xi_T)+D_T$ or the jump point $L$ on the second part, 0 on the third part.
			
				If $\xi_T<U^{'}(L-D_T)/\lambda_1$ we have that $I(\lambda_{1}\xi_T)+D_T>L$, and the Lagrangian increases continuously from $D_T$ to $I(\lambda_{1}\xi_T)+D_T$ except possibly at $L$ and decreases from then on. If $\xi_T\geq U^{'}(L-D_T)/(\lambda_{1}-\lambda_2)$, we have that $I((\lambda_{1}-\lambda_2)\xi_T)+D_T<L$. The Lagrangian increases continuously from $D_T$ to $I((\lambda_{1}-\lambda_2)\xi_T)+D_T$ and decreases from then on except possibly at the jump point $L$.
				
					By Lemma \ref{tang} we know that on  $\xi_T<U^{'}(L-D_T)/\lambda_1$, $$\phi(I(\lambda_{1}\xi_T)+D_T)-\phi(L)
		=U(I(\lambda_{1}\xi_T))-\lambda_{1}\xi_TI(\lambda_{1}\xi_T)+\lambda_{1}\xi_T(L-D_T)-U(L-D_T)>0,$$
		and  on $\xi_T\geq U^{'}(L-D_T)/(\lambda_{1}-\lambda_2)$ we have that  
		\begin{align*}
		&\phi(I((\lambda_{1}-\lambda_2)\xi_T)+D_T)-\phi(L)\\
		=&U(I((\lambda_{1}-\lambda_2)\xi_T))-(\lambda_{1}-\lambda_2)\xi_TI((\lambda_{1}-\lambda_2)\xi_T)-(\lambda_{1}-\lambda_2)\xi_TD_T\\
		&+(\lambda_{1}-\lambda_2)\xi_TL-U(L-D_T)>0.
		\end{align*}
		
		If $U^{'}(L-D_T)/\lambda_1\leq\xi_T<U^{'}(L-D_T)/(\lambda_{1}-\lambda_2)$, we have that $I(\lambda_{1}\xi_T)+D_T\leq L<I((\lambda_{1}-\lambda_2)\xi_T)+D_T$, and the Lagrangian increases from $D_T$ to the jump point $L$ and decreases from then on. Hence, $L$ is the maximizer.
		
			Therefore, the maximizer $X_T^{'}$ of the Lagrangian defined on $X_T>D_T$ is
		\begin{align}\label{lel1}
		    X_T^{'}&=(I(\lambda_{1}\xi_T)+D_T)\mathbbm{1}_{\xi_T<U^{'}(L-D_T)/\lambda_1}+L\mathbbm{1}_{U^{'}(L-D_T)/\lambda_1\leq\xi_T<U^{'}(L-D_T)/(\lambda_{1}-\lambda_2)}\nonumber\\
		    +&(I((\lambda_{1}-\lambda_2)\xi_T)+D_T)\mathbbm{1}_{\xi_T\geq U^{'}(L-D_T)/(\lambda_{1}-\lambda_2)}.
		\end{align}
			Next, we compare the maximum given by (\ref{lel1}) with $\phi(0)$ in different regions to determine the global maximum.\\
		On $\xi_T<U^{'}(L-D_T)/\lambda_1$, we have that
		\begin{align}\label{13}
		\phi(I(\lambda_{1}\xi_T)+D_T)-\phi(0)&=U(I(\lambda_{1}\xi_T))-\lambda_{1}\xi_TI(\lambda_{1}\xi_T)-\lambda_{1}\xi_TD_T+\lambda_2\xi_TL\nonumber\\
		&=\max\limits_{X_T} \bigg\{U(X_T-D_T)-X_T\lambda_{1}\xi_T+\lambda_{1}\xi_T\frac{\lambda_2}{\lambda_{1}}L\bigg\}.
		\end{align}
		Lemma \ref{tang} provides the necessary and sufficient condition for the existence of the zero root of (\ref{13}) denoted by $\xi_a$. We have that $I(\lambda_{1}\xi_T)+D_T$ is the global maximizer if  $\omega \in \{\xi_T<U^{'}(L-D_T)/\lambda_1\} \bigcap \{\xi_T<\xi_a\}$.
		
			Similarly on $U^{'}(L-D_T)/\lambda_1\leq\xi_T<U^{'}(L-D_T)/(\lambda_{1}-\lambda_2)$, we have that $\phi(L)-\phi(0)=U(L-D_T)-\lambda_{1}\xi_TL+\lambda_2\xi_TL$ and $\xi_b=U(L-D_T)/(\lambda_{1}-\lambda_2)L$ is the zero root. Therefore, if $\omega \in \{U^{'}(L-D_T)/\lambda_1<\xi_T<U^{'}(L-D_T)/(\lambda_{1}-\lambda_2)\}\bigcap\{\xi_T<\xi_b\}$, $L$ is the global maximizer.
			
				Further, on $\xi_T\geq U^{'}(L-D_T)/(\lambda_{1}-\lambda_2)$, we get that $\phi(I((\lambda_{1}-\lambda_2)\xi_T)+D_T)-\phi(0)
		=U(I((\lambda_{1}-\lambda_2)\xi_T))-(\lambda_{1}-\lambda_2)\xi_TI((\lambda_{1}-\lambda_2)\xi_T)-(\lambda_{1}-\lambda_2)\xi_TD_T$ of which the zero root is $\xi_c:=U^{'}(\wh D_T-D_T)/(\lambda_{1}-\lambda_2)$. Hence, if $\omega \in \{\xi_T>U^{'}(L-D_T)/(\lambda_{1}-\lambda_2)\}\bigcap\{\xi_T<\xi_c\}$, $I((\lambda_{1}-\lambda_2)\xi_T)+D_T$ is the global maximizer. In any other case, zero is the global maximizer.
		
			For the sake of clarity, we formulate a table of all the critical values of $\xi_T$.
		\begin{center}
			\begin{tabular}{l|l|l}
				\hline
				$U^{'}(L-D_T)/\lambda_{1}$&$U^{'}(L-D_T)/(\lambda_{1}-\lambda_2)$&$U^{'}(\wh D_T-D_T)/\lambda_{1}$\\
				\hline
				$\xi_a:=$&$\xi_b:=$&$\xi_c:=$\\
				
				zero root of (\ref{13})&$U(L-D_T)/L(\lambda_{1}-\lambda_2)$&$U^{'}(\wh D_T-D_T)/(\lambda_{1}-\lambda_2)$\\
				\hline
					\end{tabular}
		\end{center}
		The utility function is concave on $x>D_T$. Assuming $L>\wh D_T$ we have that $\xi_a>U^{'}(\wh D_T-D_T)/\lambda_1>U^{'}(L-D_T)/\lambda_1$, $\xi_c>U^{'}(L-D_T)/(\lambda_{1}-\lambda_2)$ and $\xi_b>U^{'}(L-D_T)/(\lambda_{1}-\lambda_2)$ . This implies that $\{\xi_T<U^{'}(\wh D_T-D_T)/\lambda_1\} \bigcap \{\xi_T<\xi_a\}=\{\xi_T<U^{'}(\wh D_T-D_T)/\lambda_1\}$, $\{U^{'}(\wh D_T-D_T)/\lambda_1<\xi_T<U^{'}(L-D_T)/(\lambda_{1}-\lambda_2)\}\bigcap\{\xi_T<\xi_b\}=\{U^{'}(\wh D_T-D_T)/\lambda_1<\xi_T<U^{'}(L-D_T)/(\lambda_{1}-\lambda_2)\}$, 
		and $\{\xi_T>U^{'}(L-D_T)/(\lambda_{1}-\lambda_2)\}\bigcap\{\xi_T<U^{'}(\wh D_T-D_T)/\lambda_1\}=\{U^{'}(\wh D_T-D_T)/\lambda_1>\xi_T>U^{'}(L-D_T)/(\lambda_{1}-\lambda_2)\}$. Therefore, 
			\begin{align}\label{17}
			    X_T^{ESQ}&=(I(\lambda_{1}\xi_T)+D_T)\mathbbm{1}_{\xi_T<U^{'}(\wh D_T-D_T)/\lambda_1}+L\mathbbm{1}_{U^{'}(\wh D_T-D_T)/\lambda_1\leq\xi_T<U^{'}(L-D_T)/(\lambda_{1}-\lambda_2)}\nonumber\\
			    +&(I((\lambda_{1}-\lambda_2)\xi_T)+D_T)\mathbbm{1}_{U^{'}(L-D_T)/(\lambda_{1}-\lambda_2)\leq\xi_T<\xi_c},
			\end{align}
					is the global maximizer of the Lagrangian in this case.
		
		If $L\leq\wh D_T$, we have that $\xi_c<U^{'}(L-D_T)/(\lambda_{1}-\lambda_2)$. Hence, $\{\xi_T>U^{'}(L-D_T)/(\lambda_{1}-\lambda_2)\}\bigcap\{\xi_T<\xi_c\}=\emptyset$, and then $I((\lambda_{1}-\lambda_2)\xi_T)+D_T$ cannot be the maximizer.
		
			In addition, we know that $\xi_a>\xi_b$ holds because
		\begin{align*}
		0=&U(L-D_T)-\lambda_{1}\xi_bL+\lambda_2\xi_bL\\
		=&\max\limits_{X_T} \bigg\{U(X_T-D_T)-X_T\lambda_{1}\xi_a+\lambda_{1}\xi_a\frac{\lambda_2}{\lambda_{1}}L\bigg\}>U(L-D_T)-\lambda_{1}\xi_aL+\lambda_2\xi_aL.
		\end{align*}
		
			Denote by $\xi_Q$ the zero root of the function $\bbe[\xi_T L\mathbbm{1}_{\xi_T>\xi_Q}=\epsilon_\bbq$.
		If $\xi_{Q}<U^{'}(L-D_T)/\lambda_1$, $\lambda_2$ is chosen such that $\xi_a=\xi_Q$. If $\xi_{Q}\geq U^{'}(\wh D_T-D_T)/\lambda_1$, $\lambda_2$ is chosen such that $\xi_b=\xi_{Q}$. We have proved that (\ref{35}) and (\ref{36}) are the global maximizer of the Lagrangian, respectively.
		
		Now we show that Proposition \ref{solu:esq} is the optimal solution. Suppose there exists another feasible solution $Y_T$. Consider
			\begin{align}
		&\bbe\left[U\left((X_T^{ESQ}-D_T)^+\right)\right]-\bbe\left[U\left((Y_T-D_T)^+\right)\right]\nonumber\\
		=&\bbe\left[\left(U(X_T^{ESQ}-D_T)^+\right)\right]-\bbe\left[U\left((Y_T-D_T)^+\right)\right]-\lambda_{1} x_0+\lambda_{1} x_0-\lambda_2\epsilon+\lambda_2\epsilon\nonumber\\
		\geq& \bbe\left[\left(U(X_T^{ESQ}-D_T)^+\right)\right]-\bbe\left[U\left((Y_T-D_T)^+\right)\right]-\bbe[\lambda_{1}\xi_TX_T^{ESQ}]+\bbe[\lambda_{1}\xi_TY_T]\label{equ28}\\
		&-\bbe[\lambda_2\xi_T(L-X_T^{ESQ})\mathbbm{1}_{X_T^{ESQ}<L}]+\bbe[\lambda_2\xi_T(L-Y_T)\mathbbm{1}_{Y_T<L}]\nonumber\\
		=&\bbe\left[\max\limits_{X_T}\left\{U\left((X_T-D_T)^+\right)-\lambda_{1}\xi_TX_T-\lambda_2\xi_T(L-X_T)\mathbbm{1}_{X_T<L}\right\}\right]\nonumber\\
		&-\bbe\left[U\left((Y_T-D_T)^+\right)-\lambda_{1}\xi_TY_T-\lambda_2\xi_T(L-Y_T)\mathbbm{1}_{Y_T<L}\right]\geq 0.\label{equ29}	
		\end{align}
		The inequality (\ref{equ28}) holds because the optimal solution  satisfies both constraints with equality. The last inequality (\ref{equ29}) is due to the fact that the optimal solution is the argmax of the static Lagrangian. The proof is complete.\endproof
		\subsection{The optimal solution for Problem VaR}
		The optimal wealth for Problem VaR is given in the following proposition.
		\begin{proposition}\label{solu var}\normalfont 
	\begin{enumerate}[i)]
		\item 	If $L\leq D_T$, the optimal wealth  is:
			\begin{equation}
		X_T^{VaR}=\left((U^{'})^{-1}(\lambda_{1}\xi)+D_T\right)\mathbbm{1}_{\xi<\xi_{2}} \quad\text{if $\xi_{2}\leq \xi_1$,}\label{26}
			\end{equation}
		\begin{equation}
		 X^{VaR}=\left((U^{'})^{-1}(\lambda_{1}\xi)+D_T\right)\mathbbm{1}_{\xi<\xi_1}+L\mathbbm{1}_{\xi_1\leq\xi<\xi_2} \quad\text{if $\xi_1<\xi_2$,}\label{esq25}
		\end{equation}	
				
					\noindent where $\xi_1=U^{'}(\wt L-(D_T-L))/\lambda_{1}$, $\wt L$ is the tangent point of $U((X_T-D_T+L)^+)$, $\xi_2$ is defined through $\bbp(\xi>\xi_2)=\alpha$, 
				and $\lambda_1$ is obtained by solving $\bbe[\xi X_T^{VaR}]=x_0$.
					\item If $D_T<L\leq\wh D_T$, the optimal wealth is:
						\begin{align}
				X_T^{VaR}=&\left((U^{'})^{-1}(\lambda_{1}\xi)+D_T\right)\mathbbm{1}_{\xi<\xi_{2}}\quad\text{if $\xi_{2}<\xi_1$,}\label{35}\\
				X_T^{VaR}=&\left((U^{'})^{-1}(\lambda_{1}\xi)+D_T\right)\mathbbm{1}_{\xi<\xi_1}+L\mathbbm{1}_{\xi_1\leq\xi<\xi_2}\quad\text{if $\xi_1\leq\xi_2$,}\label{36}
				\end{align}
					where $\xi_1=U^{'}(L-D_T)/\lambda_{1}$ and $\xi_2$ is defined through $\bbp(\xi>\xi_2)=\alpha$, 
				and $\lambda_1$ is obtained by solving $\bbe[\xi X_T^{VaR}]=x_0$.
		\item If $L>\wh D_T$, the optimal wealth is:
		\begin{footnotesize}
		             \begin{equation}
		X_T^{VaR}=\left((U^{'})^{-1}(\lambda_{1}\xi)+D_T\right)\mathbbm{1}_{\xi<\xi_1}+L\mathbbm{1}_{\xi_1\leq\xi<\xi_2}+\left((U^{'})^{-1}(\lambda_{1}\xi)+D_T\right)\mathbbm{1}_{\xi_2<\xi<\xi_3}\quad\text{if $\xi_{2}<\xi_3$,}\label{35}
		\end{equation}
			\begin{equation}	X_T^{VaR}=\left((U^{'})^{-1}(\lambda_{1}\xi)+D_T\right)\mathbbm{1}_{\xi<\xi_1}+L\mathbbm{1}_{\xi_1\leq\xi<\xi_2}\quad\text{if $\xi_3\leq\xi_2$,}\label{36}
				\end{equation}
					\end{footnotesize}
					\noindent where $\xi_1=\frac{ U^{'}(L-D_T)}{\lambda_1}$, $\xi_2$ is defined through $\bbp(\xi>\xi_2)=\alpha$,
			$\xi_3=\frac{U^{'}(\wh D_T-D_T)}{\lambda_1}$, $\lambda_1$ is obtained by solving the equations $\bbe[\xi X_T^{VaR}]=x_0$.
	\end{enumerate}
\end{proposition}
\proof
With a similar argument as in the proof for the optimal solution for Problem ESQ, and the help of Lemma \ref{tang}, we can show that Proposition \ref{solu var} gives the optimal solution for Problem VaR.\endproof
	\subsection{Proof of Proposition \ref{proposition:AVaR}}
	The proof of Proposition \ref{proposition:esp} (the optimal solution for Problem ESP) is similar to the proof in Section \ref{sec:proofesq} (the optimal solution for Problem ESQ). In this section, we provide the proof for Proposition \ref{proposition:AVaR}, i.e., the optimal solution for Problem AVaR is equivalent to the optimal solution for Problem ESP.

\proof
Let $U^{max}$ be the supremum of the AVaR constrained problem, that is 
$$U^{max}:=\sup\limits_{X_T}\bbe[U(X_T-D_T)^+],$$
where $X_T$ satisfies $\bbe[X_T\xi_T]\leq x_0$ and $\frac{1}{\alpha}\int_0^{\alpha}VaR_{X_T}(\beta)\d \beta\geq L^{AVaR}$. Note that for each given AVaR constraint, there is no feasible solution if the initial wealth is too small. On the other hand, if the initial wealth is too large, the unconstrained benchmark solution is the optimal one. We do not discuss these trivial cases here but only focus on the case where there is a non-trivial solution. 

Because of Lemma \ref{lemma:AVaR}, $U^{max}$ is also the supremum of the ES-constrained problem, that is
\begin{align*}
U^{max}&=\sup\limits_{L\geq L^{AVaR}}\bbe[U(X_T^{ES,L}-D_T)^+],
\end{align*}
where $X_T^{ES,L}$ is the optimal solution for Problem ESP with the constraint $\bbe[(L-X_T^{ES,L})^+]=(L-L^{AVaR})\alpha$.

We first show that the set of the implicit threshold is close and bounded for each given initial wealth $x_0$. Then, we show that the supremum of the expected utility is attainable at the implicit threshold $L^{*}$ since the expected utility is continuous in the implicit threshold.

 Because of Lemma \ref{lemma:AVaR}, for each given initial wealth $x_0$, it is sufficient to consider the implicit threshold $L^{*}$ such that there exists an attainable wealth $X_T^{*}$ with $L^{*}$ being its $\alpha-$quantile, i.e., $L^{*}=q_{X_T^{*}}(\alpha)$. Thus an upper bound of $L^{*}$ is given by 
$$\max\limits_{X_T}q_{X_T}(\alpha),\quad\text{subject to}\quad \bbe[X_T\xi]\leq x_0.$$
$$$$
The solution is given by
$$X_T=q_{X_T}(\alpha)\mathbbm{1}_{\xi\leq q_{\xi}(1-\alpha)},$$
where $q_{\xi}(1-\alpha)$ denotes the $(1-\alpha)$-quantile of the state price density. 

Thus, in a complete financial market, we can use the initial wealth to calculate explicitly the upper bound of the $\alpha-$quantile of a feasible solution, that is 
$$\wh q_{X_T}(\alpha)=x_0/\bbe[\xi\mathbbm{1}_{\xi\leq q_{\xi}(1-\alpha)}],$$ 
which is also the upper bound of the implicit threshold $L$, i.e.,
$\wh L=\wh q_{X_T}(\alpha)$
since $L\equiv q_{X_T}(\alpha)$ holds in the equivalent problem. Hence, the supremum utility $U^{max}$ can be rewritten as 
$$U^{max}=\sup\limits_{L^{AVaR}\leq L\leq\wh L}\bbe[U(X_T^{ES,L}-D_T)^+].$$

From the explicit solutions of the ESP constrained problem we discussed in the previous section, we know that the function: $L\rightarrow\sup\bbe[(U^{ES,L}-D_T)^+]$ is continuous since $L\rightarrow X_T^{ES,L}$ is continuous a.s. Therefore, we conclude that there exists $L^{*}$ such that the supremum utility $U^{max}$ is attained. By Lemma \ref{lemma:AVaR}, $U^{max}$ is also the maximum utility for a feasible solution under the AVaR constraint, which implies that the optimal solution under the AVaR constrained coincides with the ES-constrained solution, that is, $X_T^{AVaR}=X_T^{ES,L^{*}}$.\endproof
\section{The existence of the Lagrangian multipliers}
Finally, we show the existence of the Lagrangian multipliers. In the previous section, we applied the Lagrangian approach assuming that the Lagrangian multipliers exist and are non-negative. In this section, we show that the Lagrangian multipliers truly exist.

		\begin{lemma}\label{lambda}
			\normalfont
			In the constrained optimization problem with any given feasible initial wealth ($x_0\geq x_0^{\min}$), the corresponding Lagrangian multipliers exist such that the budget constraint and the risk constraint are simultaneously binding.
		\end{lemma}
		We classify all the optimal solutions into two groups. The first group contains solutions which depend only on the budget Lagrangian multiplier $\lambda_1$ and the second group contains solutions which depend on two Lagrangian multipliers $\lambda_1$ and $\lambda_2$ (see Propositions \ref{proposition:esp}, \ref{solu:esq} and \ref{solu var}). We see that all the optimal solutions belong to the first group except for the solutions for Problems ESP and ESQ constraints in the case $L>\wh D_T$. In the first group, $\lambda_2$ is a function of $\lambda_1$ assuming the risk constraint is binding. Therefore, as long as we show that $\lambda_1$ exists, $\lambda_2$ exists as well.
		
			\proof We first prove that the budget Lagrangian multiplier $\lambda_1$ in the first group of solutions exists. Take the two-region solution for Problem ESQ $X^{ESQ}(\lambda_{1},\xi)=(U^{'})^{-1}((\lambda_{1}\xi)+D_T)\mathbbm{1}_{\xi<\xi_{2}}$ in the case $L\leq D_T$ as an example. Recall that $\lambda_1$ is the solution of $\bbe[\xi X_T^{ESQ}]=x_0$ and $\xi_2$ is the solution of $\bbe[\xi L\mathbbm{1}_{\xi>\xi_2}]=\epsilon_{\bbq}$.
		
		We first show that the map $\varphi : \lambda \rightarrowtail \bbe[\xi_TX^{ESQ}(\lambda,\xi)]$ is a strictly decreasing, continuous and surjective function from $(0,\infty)$ to $(x_0^{\min},\infty)$.
		
		Showing that $\varphi$ is a strictly decreasing function is equivalent to show that for all $\lambda_a>\lambda_b>0, \varphi(\lambda_a)<\varphi(\lambda_b)$. Define $l(\lambda):=(I(\lambda\xi_T)+D_T)\mathbbm{1}_{\xi<\xi_{2}}$. Then we have $\varphi(\lambda_a)=\bbe[\xi l(\lambda_a)]$ and $\varphi(\lambda_b)=\bbe[\xi l(\lambda_b)]$. We know that  $I(\lambda\xi)+D_T$ is a strictly decreasing function in $\lambda$ almost surely. Hence, we conclude that $l(\lambda_a) \leq l(\lambda_b)$ and $\bbp(l(\lambda_a)<l(\lambda_b))>0$ as long as $\{\xi<\xi_{2}\}$ is not an empty set. Further, we have $ \varphi(\lambda_a)=\bbe[\xi_Tl(\lambda_a)]<\bbe[\xi_Tl(\lambda_b)]=\varphi(\lambda_b)$. Thus, $\varphi$ is a strictly decreasing function.
		
			It is intuitive to see that $l(\cdot)$ is a continuous function except for countable many points. Therefore, $\varphi$ is a continuous function in $\lambda$ a.s.
		
		In addition, if $\lambda$ tends to $\infty$ then $\varphi(\lambda)$ tends to $x_0^{\min}$ and if $\lambda$ tends to zero then $\varphi$ tends to $\infty$. Therefore, $\varphi$ is a strictly decreasing, continuous and surjective function from $(0,\infty)$ to $(x_0^{\min},\infty)$. Hence, for each fixed $x_0\geq x_0^{\min}$, the equation $\bbe[\xi X_T^{ESQ}]=x_0$ has a unique zero root which is $\lambda_1$.
		
		Letting $I(\cdot)$ denote the inverse function of $(U^{'})^{-1}$, from equation (\ref{esq1}), we know that $$\lambda_2=\frac{-(U(I(\lambda_{1}\xi_{2}))-I(\lambda_{2}\xi_{2})\lambda_{1}\xi_{2}-\lambda_{1}\xi_{2}D_T)}{\xi_{2}L}$$ and the ESQ-constraint is binding. 
		
		The other cases in the first group can be proved with similar arguments. 
		
		Next we prove that the two Lagrangian multipliers exist simultaneously such that the budget constraint and the ESQ constraint are binding simultaneously in the case $L>\wh D_T$.
		
			In the first step of the proof we show that for a fixed $\lambda_{1}$ and for $0<\lambda_2 \leq \lambda_{1}$, the second constraint always holds.
		
		If $\lambda_2$ tends to $\lambda_{1}$, $ U^{'}(L-D_T)/(\lambda_{1}-\lambda_2)$ converges to $\infty$. Hence, the optimal solution converges to
		\begin{equation*}
		X^{ESQ}=\begin{cases}
		I(\lambda_{1}\xi)+D&\quad\text{if}\quad\xi<U^{'}(L-D_T)/\lambda_{1},\\
		L&\quad\text{if}\quad U^{'}(L-D_T)/\lambda_{1}\leq\xi.\\
		\end{cases}
		\end{equation*}
		This implies that $\bbe[(L-X_T^{ESQ})\xi\mathbbm{1}_{X_T^{ESQ}<L}]-\epsilon_\bbq=-\epsilon_\bbq<0$, which obviously satisfies the ESQ constraint. 
		
		If $\lambda_2=0$, the optimal solution converges to
		\begin{equation*}
		X^{ES_\bbq}=(I(\lambda_{1}\xi)+D_T)\mathbbm{1}_{\xi<U^{'}(\wh D_T-D_T)/\lambda_{1}},
		\end{equation*}
		which is the benchmark solution.
			For a given $\lambda_{1}$, $\bbe[(L-X_T^{ESQ})\xi\mathbbm{1}_{X_T^{ESQ}<L}]-\epsilon_\bbq$ is a continuous and 	 decreasing function in $\lambda_2$ and thus bijective. By the intermediate value theorem, the zero root of the function exists. Let $\lambda_2(\lambda_{1})$ denote the zero root, which is a function of $\lambda_{1}$.
		
		Denoting by $X_T^{ESQ,\lambda_{1},\lambda_2(\lambda_{1})}$ the optimal terminal wealth, we have that
		\begin{align*}
		\bbe[X_T^{ESQ,\lambda_{1},\lambda_2(\lambda_{1})}\xi_T]=&\bbe[X_T^{ESQ,\lambda_{1},\lambda_2(\lambda_{1})}\xi\mathbbm{1}_{X_T^{ESQ}<L}]+\bbe[X_T^{ESQ,\lambda_{1},\lambda_2(\lambda_{1})}\xi\mathbbm{1}_{X_T^{ESQ,\lambda_{1},\lambda_2(\lambda_{1})}\geq L}]\\
		=&L\bbe[\xi\mathbbm{1}_{X_T^{ESQ,\lambda_{1},\lambda_2(\lambda_{1})}<L}]-\epsilon_\bbq+\bbe[X_T^{ESQ,\lambda_{1},\lambda_2(\lambda_{1})}\xi\mathbbm{1}_{X_T^{ESQ,\lambda_{1},\lambda_2(\lambda_{1})}> L}]\\
		&+\bbe[X_T^{ESQ}\xi\mathbbm{1}_{X_T^{ESQ,\lambda_{1},\lambda_2(\lambda_{1})}= L}]\\
		=&L\bbe[\xi\mathbbm{1}_{X_T^{ESQ,\lambda_{1},\lambda_2(\lambda_{1})}<L}]-\epsilon_\bbq+\bbe[X_T^{ESQ,\lambda_{1},\lambda_2(\lambda_{1})}\xi\mathbbm{1}_{X_T^{ESQ,\lambda_{1},\lambda_2(\lambda_{1})}>L}]\\
			&+\bbe[L\xi\mathbbm{1}_{X_T^{ESQ,\lambda_{1},\lambda_2(\lambda_{1})}= L}]\\
		=&-\epsilon_\bbq+\bbe[L\xi\mathbbm{1}_{X_T^{ESQ,\lambda_{1},\lambda_2(\lambda_{1})}\leq L}]+\bbe[X_T^{ESQ,\lambda_{1},\lambda_2(\lambda_{1})}\xi\mathbbm{1}_{X_T^{ESQ,\lambda_{1},\lambda_2(\lambda_{1})}>L}]\\
		=&L\bbe[\xi]+\bbe[(X_T^{ESQ,\lambda_{1},\lambda_2(\lambda_{1})}-L)\xi\mathbbm{1}_{X_T^{ESQ,\lambda_{1},\lambda_2(\lambda_{1})}>L}]-\epsilon_\bbq\\
		=&Le^{-rT}-\epsilon_\bbq+\bbe[(X_T^{ESQ,\lambda_{1},\lambda_2(\lambda_{1})}-L)\xi\mathbbm{1}_{X_T^{ESQ,\lambda_{1},\lambda_2(\lambda_{1})}>L}].
		\end{align*}
		By equation (\ref{17}) we have  $\mathbbm{1}_{X_T^{ESQ,\lambda_{1},\lambda_2(\lambda_{1})}>L}=\mathbbm{1}_{\xi_T<U^{'}(L-D_T)/\lambda_{1}}$. If $\lambda_{1}$ tends to zero, $X_T^{ESQ,\lambda_{1},\lambda_2(\lambda_{1})}$ converges to $\infty$ and $\bbe[X_T^{ESQ}\xi]$ converges to $\infty$. On the other hand if $\lambda_{1}$ tends to $\infty$, $\bbe[X_T^{ESQ}\xi]$ converges to $Le^{-rT}-\epsilon_\bbq$.
			Moreover, $\bbe[X_T^{ESQ}\xi]$ is a continuous function in $\lambda_{1}$ and thus as long as $x_0>Le^{-rT}-\epsilon_\bbq$, by the intermediate value theorem there exists $\lambda_{1}$ such that $\bbe[X_T^{ESQ,\lambda_{1},\lambda_2(\lambda_{1})}\xi]=x_0$.
		
		Following a similar procedure, we can show that $\lambda_{1}$ and $\lambda_2$ exist in the optimal solution for Problem ESP in the case $L>\wh D_T$. Recall that the optimal solution is given by
		\begin{equation}\label{17}
				X_T^{ESP}=\left((U^{'})^{-1}(\lambda_1\xi)+D_T\right)\mathbbm{1}_{\xi<\xi_1}+L\mathbbm{1}_{\xi_1\leq\xi<\xi_2}+\left((U^{'})^{-1}(\lambda_1\xi-\lambda_2)+D_T\right)\mathbbm{1}_{\xi_2\leq\xi<\xi_3},
				\end{equation}
				
					\noindent where $\xi_1=\frac{ U^{'}(L-D_T)}{\lambda_1}$, $\xi_2=\frac{ U^{'}(L-D_T)+\lambda_2}{\lambda_1}$,
			$\xi_3=\frac{U^{'}(\wh D_T-D_T)+\lambda_2}{\lambda_1}$, $\lambda_1$ and $\lambda_2$ are obtained by solving the equations $\bbe[\xi X_T^{ESP}]=x_0$ and $\bbe[(L-X_T^{ESP})^+]=\epsilon_{\bbp}$.
		
		The first step is to show that given a $\lambda_{1}$, there exists a $\lambda_2$ such that the $ESP$ constraint is fulfilled, i.e., there exists a $\lambda_2$ such that $\bbe[(L-X_T^{ESP})^+]=\epsilon_\bbp$. 
		
		Note that $\lambda_2\to\infty\Rightarrow\xi_2$ and $\xi_3\to\infty$. Thus, the solution tends to 
		$$X_T^{ESP}=(I(\lambda_{1}\xi)+D_T)\mathbbm{1}_{\xi<U^{'}(L-D_T)/\lambda_1}+L\mathbbm{1}_{\xi_T\geq U^{'}(L-D_T)/\lambda_1}.$$
		
		Obviously, $\bbe[(L-X_T^{ESP})^+]=0<\epsilon_\bbp$.
		
			On the other hand, the solution tends to the benchmark solution if $\lambda_2\to 0$. In this case, we know that $\bbe[(L-X_T^{ESP})^+]>\epsilon_\bbp$. Otherwise, the constraint is redundant. Thus, by the intermediate value theorem, there exists a $\lambda_2\geq 0$ such that $\bbe[(L-X_T^{ESP})^+]=\epsilon_\bbp$. 
			\noindent Next, we show that there exists $\lambda_{1}$ such that $\bbe[X_T^{ESP}\xi]=x_0$. Note that $\lambda_{1}\to \infty$ implies that $U^{'}(L-D_T)/\lambda_1\to 0$, and the solution tends to
		$$X_T^{ESP}=L\mathbbm{1}_{\xi<\xi_2}+(I(\lambda_{1}\xi)+D_T)\mathbbm{1}_{\xi_2\leq\xi<\xi_3}.$$ 
		In this case, $\bbe[X_T^{ESP}\xi]\to x_0^{\min}$. On the other hand, $\lambda_{1}\to 0\Rightarrow\bbe[X_T^{ESP}\xi]\to\infty$. Again, by the intermediate value theorem, there exists a $\lambda_{1}$ such that $\bbe[X_T^{ESP}\xi]=x_0$.
		\endproof
\end{document}